\documentclass{article}
\usepackage{amssymb}
\usepackage{graphics}
\input xy
\xyoption{all}
\begin{document}
\title{Howe type duality for metaplectic group acting on symplectic spinor valued forms}
\author{Svatopluk Kr\'ysl \footnote{{\it E-mail address}: krysl@karlin.mff.cuni.cz}\\ {\it \small  Charles University, Sokolovsk\'a 83, Praha 8, Czech Republic} \\{\it \small and
Humboldt-Universit\"{a}t zu Berlin, Unter den Linden 6, Berlin, Germany.}
\thanks{I am very grateful to Roger Howe for explaining me a general framework of the studied type of duality. The author of this article was supported by the grant
GA\v{C}R 201/06/P223 of the Grant Agency of the Czech Republic for
young researchers. The work is a part of the research project MSM0021620839 financed by M\v{S}MT \v{C}R.  Also supported by SFB 1096 of the DFG.}}

\maketitle

\noindent\centerline{\large\bf Abstract} Let $\lambda: \tilde{G}\to
G$ be the non-trivial double covering of the symplectic group
$G=Sp(\mathbb{V},\omega)$ of the symplectic vector space
$(\mathbb{V},\omega)$ by the metaplectic group
$\tilde{G}=Mp(\mathbb{V},\omega).$ In this case, $\lambda$ is also a
representation of $\tilde{G}$ on the vector space $\mathbb{V}$ and
thus, it gives rise to the representation of $\tilde{G}$ on the
space of exterior forms $\bigwedge^{\bullet}\mathbb{V}^*$   by
taking wedge products. Let ${\bf S}$ be the minimal globalization of
the Harish-Chandra module of the complex Segal-Shale-Weil
representation of the metaplectic group $\tilde{G}.$ We prove that
the associative commutant algebra
$\hbox{End}_{\tilde{G}}(\bigwedge^{\bullet}\mathbb{V}^*\otimes {\bf
S})$ of the metaplectic group $\tilde{G}$ acting on the ${\bf S}$-valued
exterior forms  is generated by certain representation of the super
ortho-symplectic Lie algebra $\mathfrak{osp}(1|2)$ and two distinguished operators. This establishes
a Howe type duality between the metaplectic group and the super Lie
algebra $\mathfrak{osp}(1|2).$ Also the space $\bigwedge^{\bullet}\mathbb{V}^*\otimes {\bf S}$ is decomposed wr. to the joint action of  $Mp(\mathbb{V},\omega)$ and $\mathfrak{osp}(1|2).$

{\it Math. Subj. Class.:} 22E46, 22E47, 22E45, 70G45, 81S10


\section{Introduction}

 The Howe type duality is a  result in representation theory, which enables one to see some theorems in mathematics as its special cases or consequences.
In this way, it gives a unification view into several mathematical
disciplines. This duality could be applied to derive a version of
the de Rham theorem, theorems in Hodge theory of K\"{a}hler
manifolds or the  theory of spherical harmonics. It could also be applied
in other interesting instances, like combinatorial structure of
Schubert cell decomposition of generalized flag manifolds, which
are important Klein models of many geometries, like conformal, CR,
almost Grassmannian, projective, contact projective etc. For a comprehensive treatment on the Howe type duality and
its applications, see Goodman, Wallach \cite{GW}, Goodman \cite{G},
Howe \cite{Howe} and also Rubentahler \cite{Rubenthaler}.  

We shall introduce this duality briefly now.
Let $ \rho: G \to \hbox{Aut}(W)$ be a finite dimensional representation of a reductive Lie group $G.$ Suppose $W$ is completely reducible. We want to decompose $W$ into isotypic components $W=\bigoplus_{\lambda \in \hat{G}}m_{\lambda}W^{\lambda},$ where  $\hat{G}$ is an index set, $W^{\lambda}$ is an irreducible representation of $G$, and $m_{\lambda}$ is the multiplicity of the occurrence of $W^{\lambda}$ in $W.$ To be more specific, if $G$ is a simple Lie group, the index set $\hat{G}$ equals a subset of the semi-lattice of dominant analytically integral weights (wr. to the usual choices).

We can also write the decomposition of $W$ as
 $W \simeq \bigoplus_{\lambda\in \hat{G}} W_{\lambda}\otimes W^{\lambda},$ where
$W_{\lambda}$ is a vector space of dimension $m_{\lambda}$ and $G$ acts as $1\otimes \rho(g)$ for $g\in G$ at the right hand side.
If $\hbox{dim}\,W_{\lambda} >1$ for an element $\lambda \in \hat{G},$ we say that $W$ is not multiplicity-free. The basic idea of the Howe type duality is to get a control over the multiplicities  by considering a further structure, namely the commutant algebra which will be described bellow. Before doing so, let us say a few words about the introduced notion.
The vector space $W_{\lambda}$ is usually referred to as the multiplicity space. It is easy to see, that $W_{\lambda}\simeq \hbox{Hom}_{G}(W^{\lambda},W).$\footnote{As usually for two representations $(\rho, W)$ and $(\rho',W')$ of $G,$ we denote $\hbox{Hom}_{G}(W,W'):=\{T \in \hbox{Hom}(W,W')| T\rho(g)=\rho'(g)T, \mbox{\, for all \,} g\in G\}$ and similarly for Lie algebra modules.}
The multiplicity space $W_{\lambda}$ becomes a  $\hbox{End}_{G}(W)$-module by the natural action $\sigma(A)(T)= A \circ T$ for $A\in \hbox{End}_{G}(W)$ and
$T \in W_{\lambda}.$ We shall denote $\hbox{End}_{G}(W)$ by $R^G$ shortly and call it the (associative) commutant algebra of $G$ on $W$ or sometimes, the algebra of $G$-invariants.  Then the Howe type duality can be formulated as follows:

\vspace{0.5cm}

{\bf Howe type duality. }{\it The following decomposition of $W$ $$W\simeq \bigoplus_{\lambda \in \hat{G}}W_{\lambda}\otimes W^{\lambda}$$ is a multiplicity-free decomposition into irreducible summands over $R^G \times G$ via the action $\sigma \otimes \rho.$}

\vspace{0.5cm}
By multiplicity-free, we mean that  $\lambda\neq \mu$ implies
$W^{\lambda} \not\simeq W^{\mu}$ (as $G$-modules) and $W_{\lambda} \not\simeq W_{\mu}$ (as $R^G$-modules).

The well known Schur duality is the Howe type duality in which $G$ is the general linear group $GL(\mathbb{U})$ of an $n$-dimensional complex vector space $\mathbb{U},$ $W=\mathbb{U}^{\otimes k}$ is the $k$-th tensor product of $\mathbb{U},$
$\rho$ is the $k$-th tensor product of the defining representation of $G$ and the index set $\hat{G}$ is the set of all partitions of $n$ into $k$ non-increasing non-negative numbers. In this case, the commutant  algebra $R^G$ is the group algebra of the image of the symmetric group $\mathfrak{S}_k$ under the natural representation of $\mathfrak{S}_k$ on the space $W=\mathbb{U}^{\otimes k}$ given by permuting positions in the tensor product $\mathbb{U}^{\otimes k}.$

The theory of spherical harmonics is a special case of the Howe type duality for the appropriate orthogonal group of a
Euclidean space acting on the space of polynomials defined over this vector space. The commutant algebra in this case, can be defined using certain highest weight representation of the Lie algebra $\mathfrak{sl}(2,\mathbb{C}).$
For this and the preceding example, see, e.g., Goodman, Wallach \cite{GW} and Goodman \cite{G}.

In the Clifford analysis, it is very natural to study spinor valued exterior forms and define  various analogues of the Cauchy-Riemann operator
using these spinor valued exterior forms as modules over the appropriate spin groups. In this situation, again a representation of $\mathfrak{sl}(2,\mathbb{C})$ comes into play.  See, e.g., Delanghe, Sommen, Sou\v{c}ek \cite{DSS}. 

In the paper, we shall be studying a symplectic analogue of the mentioned Clifford analysis setting. However, in the symplectic case the things are complicated, mainly because the space of "spinors" is infinite dimensional.
Let us briefly describe this  situation.
Let  $(\mathbb{V},\omega)$ be a real symplectic vector space of real dimension $2l.$ Let $\tilde{G}$ be the non-trivial double covering of the symplectic group $Sp(\mathbb{V},\omega).$ This double covering is called metaplectic group, and it is denoted by $Mp(\mathbb{V},\omega).$ For further use, let us denote the Lie algebra of $G$ by $\mathfrak{g}.$ Obviously, $\mathfrak{g} =\mathfrak{sp}(\mathbb{V},\omega) \simeq \mathfrak{sp}(2l,\mathbb{R}).$
The metaplectic group $\tilde{G}$ has a distinguished infinite dimensional unitary representation, called the Segal-Shale-Weil or the oscillator representation.  It was constructed by A. Weil using the Stone-von Neumann theorem. We will be interested in an "analytical" version of this representation. We will call this version  a metaplectic representation and its underlying vector space the space of symplectic spinors. We denote this vector space by $\bf S$. 

Let us say a few words about the name spinor we used for this representation.  
Let $(\mathbb{U},B)$ be a complex even dimensional vector space together with a non-degenerate bilinear form, and $\mathfrak{so}(\mathbb{U},B)$ be the associated complex orthogonal Lie algebra. Then one can construct the (orthogonal) spinors over $\mathfrak{so}(\mathbb{U},B)$ as the exterior algebra of a maximal isotropic subspace of $(\mathbb{U},B).$
Now, let us come to the symplectic case.
One can show that the $\mathfrak{g}^{\mathbb{C}}$-module structure of the Harish-Chandra module of $\bf S$ is equivalent to a representation of $\mathfrak{g}^{\mathbb{C}}$ on the space of complex valued polynomials defined on a Lagrangian\footnote{that is, maximal isotropic subspace of $(\mathbb{V},\omega)$} subspace $\mathbb{L}$ of $(\mathbb{V},\omega),$ i.e., to
the symmetric algebra $\oplus_{i=0}^{\infty}\odot^{i}\mathbb{L}.$ Thus, the situation is symmetric in changing the symplectic form to the non-degenerate symmetric bilinear form and changing the symmetric power into the antisymmetric one.
The action of $\mathfrak{g}^{\mathbb{C}}$ on this algebra coincides with the so called Dixmier realization of the complex symplectic Lie algebra by polynomial coefficient differential operators acting on the polynomials. For more details, see, e.g., Britten, Hooper, Lemire \cite{BHL} and Kashiwara, Vergne \cite{KV}.

Now, we are in a position to introduce the representation we shall be considering. In the studied case, the Lie group under consideration is the metaplectic group $Mp(\mathbb{V},\omega)$ and the module
$W=\bigwedge^{\bullet}\mathbb{V}^*\otimes {\bf S},$ i.e., the symplectic spinor valued exterior forms.
In this case, we prove that the commutant algebra is generated by an image of the super Lie ortho-symplectic algebra $\mathfrak{osp}(1|2)$  under a super Lie algebra representation  and two projection operators $R_{\pm}.$
Thus we prove that in this case, the Howe dual partner of the metaplectic group $Mp(\mathbb{V},\omega)$ is isomorphic to the super Lie algebra $\mathfrak{osp}(1|2).$ (We are using the word Howe dual partner in an informal way.)
The crucial tool, we are using to derive this duality, is a decomposition structure of exterior forms with values in the symplectic spinors. The mentioned decomposition was derived by a use of a theorem of Britten, Hooper and Lemire in \cite{BHL}. 

A further result, contained in this text, is a decomposition of the symplectic spinor valued exterior forms with respect to the joint action of the metaplectic group and the ortho-symplectic super Lie algebra $\mathfrak{osp}(1|2).$ The decomposition is not multiplicity-free in the sense of the formulation of the Howe type duality in this Introduction, but rather two-folded. More precisely, we decompose the mentioned module into a direct sum
$\bigoplus_{j=0}^{l} ({\bf E}_{jj}^+\otimes \mathbb{G}^j )\oplus ({\bf E}_{jj}^-\otimes \mathbb{G}^j)$ where ${\bf E}_{jj}^{\pm}$ are certain irreducible 
$Mp(\mathbb{V},\omega)$-modules and $\mathbb{G}^j$ are finite dimensional irreducible $\mathfrak{osp}(1|2)$-modules. We call this decomposition two-folded because for a given $j\in \{0,\ldots, l\},$ the representation $\mathbb{G}^j$ appears twice: once multiplied by ${\bf E}_{jj}^+$ and for the second time  by ${\bf E}_{jj}^-.$
As far as we know, the mentioned results serve a new example of Howe duality and Howe dual partners. 

The motivation of this article comes from differential geometry, namely from the theory of symplectic operators of Dirac type on manifolds with a metaplectic structure.
In the same way in which one can define Dirac, twistor and Rarita-Schwinger operators in the Riemannian spin geometry
within the spinor twisted de Rham sequence, one can do a similar construction in the symplectic case. See Habermann, Habermann \cite{HH}, Habermann \cite{KH} and Kr\'{ysl} \cite{SK}.
In physics, the symplectic spinor were used, e.g., by  Green, Hull in \cite{GH} in order to describe quantum mechanics of strings in 10 dimensional spaces; by Reuter \cite{Reuter} in the realm of the theory of Dirac-K\"{a}hler fields and also by Shale in \cite{Shale}, where the symplectic spinors were used in a quantization procedure for bosonic fields.

In the second section, some basic facts on (globalization of) Harish-Chandra modules are recalled (Theorem 1, Theorem 2).
In this section, the theorem on a decomposition of the symplectic spinor valued exterior forms into irreducible summands is derived (Theorem 9).  
In the third section, a representation, denoted by $\sigma,$ of the super Lie algebra 
$\mathfrak{osp}(1|2)$ on the space of symplectic spinor valued exterior forms is introduced. In this section, we prove that $\sigma$ is a super Lie algebra representation (Theorem 11) and that it maps the super Lie algebra $\mathfrak{osp}(1|2)$ into the commutant algebra
$\hbox{End}_{Mp(\mathbb{V},\omega)}(\bigwedge^{\bullet}\mathbb{V}^*\otimes {\bf S})$ (Corollary 13). 
In the fourth section, the commutant algebra $\hbox{End}_{Mp(\mathbb{V},\omega)}(\bigwedge^{\bullet}\mathbb{V}^*\otimes {\bf S})$ is computed. We prove that the image of the representation $\sigma$  together with two additional operators $R_{\pm}$ generate enough invariants, namely the whole commutant algebra (Theorem 17). To derive this theorem, we construct the projection operators from the space of symplectic spinor valued exterior forms onto its irreducible submodules only by using the representation $\sigma$ and the two mentioned operators.
In the fifth part, we define a family of finite dimensional irreducible representations of the super Lie algebra $\mathfrak{osp}(1|2)$ and prove the mentioned two-folded Howe type duality for the studied case (Theorem 20). To derive this theorem, certain commutator relation were crucial (Lemma 19). 
In the sixth section, a known result in symplectic spin geometry is presented (Theorem 21). This theorem was proved using the relations 
between the standard generators of $\mathfrak{osp}(1|2)$ in Kr\'ysl \cite{SK}.

\section{Symplectic spinor valued exterior forms}

In this section, we shall recall   basic results about globalization of Harish-Chandra $(\mathfrak{g},K)$-modules, which we shall use in the  Sections 4 and 5.
Further, we summarize   elementary facts on the metaplectic representation. At the end, we derive
a result on a decomposition of symplectic spinor valued exterior forms into irreducible submodules over the appropriate complex symplectic Lie algebra.

\subsection{Basic facts on globalization of Harish-Chandra modules and notation}

Let $G$ be a real reductive group in the sense of Vogan, see \cite{Vogan}.
Let us remark that typical examples of these groups are Lie subgroups of the general linear group   and their finite covers. (The metaplectic group is reductive in this sense.)
Let us also fix a maximal compact subgroup $K$ of $G.$

We will consider the category $\mathcal{R}(G),$ the objects of which are locally convex Hausdorff vector spaces with a continuous
linear action of $G,$ which is admissible and  of finite length, see, e.g., Kashiwara, Schmid \cite{KS} for more details. Let us only mention that the considered Segal-Shale-Weil representation belongs to this category.
The morphisms in the category $\mathcal{R}(G)$ are linear continuous $G$-equivariant maps between the objects. Let us denote the Lie algebra of $G$ by $\mathfrak{g}.$
For a  representation
${\bf E} \in \mathcal{R}(G)$ of $G,$ we shall denote the corresponding Harish-Chandra $(\mathfrak{g},K)$-module by $E$  and when we will only  be considering its $\mathfrak{g}$-module structure, we shall use the symbol $\mathbb{E}$ for it. We shall denote the Harish-Chandra forgetful functor from $\mathcal{R}(G)$ into the category of Harish-Chandra modules $\mathcal{HC}(\mathfrak{g},K)$ by $HC.$
Let us set a further convention. Suppose ${\bf E}, {\bf F} \in \mathcal{R}(G)$ and  a morphism $A: {\bf E} \to {\bf F}$ are given. Consider the restrictions $A_{|E}: E \to F$ and $A_{|\mathbb{E}}:\mathbb{E}\to \mathbb{F}.$  Then we shall denote these restrictions simply by $A: E \to F$ and  $A: \mathbb{E} \to \mathbb{F},$ respectively.  We shall also denote the complexification $A^{\mathbb{C}}: {\bf  E}^{\mathbb{C}} \to {\bf F}^{\mathbb{C}}$ simply  by $A: {\bf E}^{\mathbb{C}} \to {\bf F}^{\mathbb{C}}.$ We hope this causes no confusion. Further, if $U$ is a complex vector space and $U'$ is a real vector space, we mean by $U\otimes U'$ the tensor product over the real numbers, whereas if both $U$ and $U'$ are complex, then by $U\otimes U'$ we mean the tensor product over the complex numbers.
Let us also stress that even if we are considering modules over a real Lie algebra, the modules are always assumed to be  complex.

We will need the following theorem  which is a special case of a theorem of Schmid.

{\bf Theorem 1:}  For a real reductive Lie group $G,$
 there exists a left adjoint   functor   to the Harish-Chandra forgetful functor $HC: \mathcal{R}(G) \to \mathcal{HC}(\mathfrak{g},K).$

{\it Proof.} See  Schmid \cite{Schmid}. $\Box$

The adjoint functor from the previous theorem is usually called minimal globalization or minimal Schmid globalization functor and will be denoted by $mg.$ In Kashiwara, Schmid \cite{KS}, the case of real points of a complex algebraic group is treated.

In the next theorem, the homomorphisms of $(\mathfrak{g},K)$-modules and
that ones  of $\mathfrak{g}$-modules are related.

{\bf Theorem 2:} Let $G$ be a real reductive Lie group, $K$ a maximal compact subgroup of $G,$ and ${\bf X}_1, {\bf X}_2 \in \mathcal{R}(G).$
If $K$ is connected, then there exists a linear bijection
$$\mbox{Hom}_\mathfrak{g}(\mathbb{X}_1,\mathbb{X}_2) \simeq \mbox{Hom}_{(\mathfrak{g},K)}(X_1,X_2).$$

{\it Proof.} See, e.g., Baldoni \cite{Baldoni}.$\Box$

We shall need the following generalization of the Schur lemma.

{\bf Theorem 3:} Let $\mathfrak{g}$ be simple complex Lie algebra and $\mathbb{X}$ be an irreducible complex $\mathfrak{g}$-module. Then every $\mathfrak{g}$-endomorphism is a scalar, i.e., $\hbox{End}_{\mathfrak{g}}(\mathbb{X}) \simeq \mathbb{C} \, \mbox{Id}_{|\mathbb{X}}.$

{\it Proof.} See  Theorem 2.6.5. and Theorem 2.6.6. in Dixmier \cite{D}. $\Box$

{\bf Remark:} From this theorem, one can easily derive that for $\mathfrak{g}$ complex simple Lie algebra, $\mathbb{X}_1$ and $\mathbb{X}_2$ irreducible $\mathfrak{g}$-modules, $\mbox{dim} \, \mbox{Hom}_{\mathfrak{g}}(\mathbb{X}_1,\mathbb{X}_2)\leq 1$ and the equality holds iff $\mathbb{X}_1\simeq \mathbb{X}_2$ as $\mathfrak{g}$-modules.

Let $\mathfrak{g}$ be a simple complex Lie algebra. Choose a Cartan subalgebra $\mathfrak{h}$ of $\mathfrak{g}$ and a set of positive roots $\Phi^+.$
We denote the irreducible highest weight complex $\mathfrak{g}$-module with a highest weight $\lambda \in \mathfrak{h}^*$ by $L(\lambda).$ If $\lambda$ happens to be dominant and integral with respect to the choice
$(\mathfrak{h}, \Phi^+),$ we denote $L(\lambda)$ by $F(\lambda)$ suggesting the fact that the module $L(\lambda)$ is finite dimensional.  For an integral dominant weight $\lambda$ wr. to
$(\mathfrak{h},\Phi^+),$ we denote the set of all weights of the
irreducible representation $F(\lambda)$ by $\Pi(\lambda).$

\subsection{Segal-Shale-Weil representation and symplectic spinors}

Consider a $2l$ dimensional real symplectic vector space
$(\mathbb{V},\omega).$ Let $\mathbb{V} = \mathbb{L}\oplus \mathbb{L}'$ be a direct sum decomposition of the vector space $\mathbb{V}$ into two Lagrangian subspaces $\mathbb{L}$ and $\mathbb{L}'$ of $(\mathbb{V},\omega).$ Let $\{e_i\}_{i=1}^{2l}$ be an adapted symplectic basis of $(\mathbb{V}=\mathbb{L}\oplus \mathbb{L}',\omega),$ i.e., $\{e_i\}_{i=1}^{2l}$ is a symplectic basis of $(\mathbb{V},\omega)$ and $\{e_i\}_{i=1}^{2l}\subseteq \mathbb{L}$ and $\{e_{i}\}_{i=l+1}^{2l}\subseteq \mathbb{L}'.$ Because the notion of a symplectic basis is not unique, let us fix one for a later use.
We call a basis $\{e_i\}_{i=1}^{2l}$ symplectic basis of $(\mathbb{V},\omega)$ if for $\omega_{ij}:=\omega(e_i,e_j),$ we have $\omega_{ij}=1$ if an only if $i\leq l$ and $j=i+l$;  $\omega_{ij}=-1$ if and only if $i>l$ and $j=i-l;$ and finally, $\omega_{ij}=0$ in other cases. For $i,j =1,\ldots, 2l,$ let us choose numbers $\omega^{ij}$ defined by $\sum_{k=1}^{2l}\omega_{ik}\omega^{kj}=\delta_i^j.$ Let us reserve the symbol $\{\check{e}_i\}_{i=1}^{2l}$ for the $\omega$-dual basis of $(\mathbb{V},\omega)$ to the basis $\{e_i\}_{i=1}^{2l},$ i.e.,
$\check{e}_i=\sum_{j=1}^l\omega^{ij}e_j$ for $i=0,\ldots, l.$
The basis of $\mathbb{V}^*$ dual to the basis $\{e_i\}_{i=1}^{2l}$
will be denoted by $\{\epsilon^i\}_{i=1}^{2l}.$ Thus, we have $\omega=\epsilon^1\wedge\epsilon^{l+1}+\ldots+\epsilon^{l}\wedge \epsilon^{2l}.$ We shall use this explicit form in the Section 4.

For the symplectic vector space $(\mathbb{V},\omega),$ we fix a symplectic group $G=Sp(\mathbb{V},\omega).$ 
Let us denote the maximal compact subgroup of $G$ by $K$. It is easy to see that $K \simeq U(l).$ Because the first homotopy group of $G$ is isomorphic to the first homotopy group of $K$ which in turn, is isomorphic to the (infinite cyclic) group $\mathbb{Z},$ there exists a non-trivial double covering space $\tilde{G}$ of $G.$
Let us denote this non-trivial double covering by $\lambda: \tilde{G}\to G.$ The group $\tilde{G}$ is usually called the  metaplectic group and the symbol $Mp(\mathbb{V},\omega)$ is used for it. In the paper, we shall find the Howe dual partner for this group and its certain infinite dimensional representation. We shall also need the $\lambda$-preimage of $K$ in $\tilde{G},$ which we denote by $\tilde{K}.$
One can easily show that
$$\tilde{K}\simeq \widetilde{U(l)}:=\{(g,z) \in U(l)\times \mathbb{C}^{\times}| \mbox{det}g=z^2\}.$$
This group is obviously connected. For more details on $\tilde{K}$, see Tirao, Vogan, Wolf \cite{TVW}.

For the Lie algebra of $\mathfrak{g}=\mathfrak{sp}(\mathbb{V},\omega)$ of $G,$
we have $\mathfrak{g}\simeq \mathfrak{sp}(2l,\mathbb{R}).$ Let us also mention that the Lie algebra of $\tilde{G}$ is  isomorphic to $\mathfrak{g}.$
For some technical reasons (uniformity of presentation of our results), let us suppose that $l\geq 3,$  and consider this choice holds throughout the rest of this article.
 We shall denote the complexification of $\mathbb{V}$  by $\mathbb{V}^{\mathbb{C}}$ and the complexification of the Lie algebra $\mathfrak{g}=\mathfrak{sp}(\mathbb{V},\omega)$ by $\mathfrak{g}^{\mathbb{C}}.$
(The complexified symplectic form on $\mathbb{V}^{\mathbb{C}}$ 
will be denoted by $\omega$.) If we choose a Cartan subalgebra
$\mathfrak{h}^{\mathbb{C}}\subseteq \mathfrak{g}^{\mathbb{C}}$ and a
system of positive roots $\Phi^+ \subseteq
(\mathfrak{h}^{\mathbb{C}})^*,$ then the set of fundamental weights
$ \{ \varpi_i \}_{i=1}^l$ is uniquely determined. This set
determines a basis $\{\epsilon_i\}_{i=1}^l$ of
$(\mathfrak{h}^{\mathbb{C}})^*$ given by
$\varpi_i=\sum_{j=1}^i\epsilon_j$ for $i=1,\ldots, l.$\footnote{Let
us mention that we are distinguishing the bases
$\{\epsilon^i\}_{i=1}^{2l} \subseteq \mathbb{V}^*$ and
$\{\epsilon_i\}_{i=1}^l \subseteq (\mathfrak{h}^{\mathbb{C}})^{*}$ from each other,
although one could identify their first $l$ members by a choice of the first
inclusion from the chain of inclusions $\mathfrak{h}^{\mathbb{C}}\subseteq \mathfrak{g}^{\mathbb{C}} \subseteq (\mathbb{V}^{*}\otimes \mathbb{V})^{\mathbb{C}}.$} It is well known that the set of roots for $(\mathfrak{g}^{\mathbb{C}},\mathfrak{h}^{\mathbb{C}})$ is given by 
$\Phi=\{\pm(\epsilon_i \pm \epsilon_j)| 1\leq i<j\leq 2l\} \cup \{\pm 2\epsilon_i|i=1,\ldots, l\}.$
We choose the set $\Phi^+=\{\epsilon_i-\epsilon_j| 1\leq i < j\leq l\} \cup\{\epsilon_i+\epsilon_j|1\leq i \leq j \leq 2l\}$ of positive roots of $(\mathfrak{g}^{\mathbb{C}},\mathfrak{h}^{\mathbb{C}}).$ The set 
 $\Delta=\{\epsilon_{i}-\epsilon_{i+1}|i=1,\ldots,l-1\}\cup\{2\epsilon_l\}$ is the set of simple roots for $(\mathfrak{g}^{\mathbb{C}},\mathfrak{h}^{\mathbb{C}},\Phi^+).$ 
For $\lambda=\sum_{i=1}^l\lambda_i\epsilon_i,$ we shall sometimes denote $L(\lambda)$ by $L(\lambda_1,\ldots,\lambda_l);$ sometimes also without the commas.
  
 For $i,j=1,\ldots, 2l,$ denote the matrix having  entry $1$ at the place $(i,j)$ an zeros at other positions by $E_{i,j}.$
Let us write the Chevalley basis elements of $\mathfrak{g}^{\mathbb{C}}=\mathfrak{sp}(\mathbb{V}^{\mathbb{C}},\omega)$ wr. to the basis 
$B:=\{e_{i}\}_{i=1}^{2l}\subseteq \mathbb{V}$ of $\mathbb{V}.$  We have then
\begin{eqnarray}
{[X_{\epsilon_i-\epsilon_j}]}_{B} &=& E_{i,j}-E_{j+l,i+l}, \quad 1 \leq i < j \leq l \\
{[X_{2\epsilon_i}]}_{B} &=& E_{i,l+i}, \quad i=1,\ldots, l \\
{[Y_{\mu}]}_B &=& {[X_{\mu}]}_B^T, \quad \mu \in \Delta\\
{[H_{\epsilon_i -\epsilon_{i+1}}]}_B &=& E_{i,i}-E_{j,j}+E_{l+j,l+j}-E_{l+i,l+i},\quad i=1,\ldots,l-1\\
{[H_{2\epsilon_l}]}_B &=& E_{l,l}-E_{2l,2l}
\end{eqnarray}
See Britten, Hooper, Lemire \cite{BHL} for more details.

Now, we shall introduce some necessary conventions we shall be using when we will be dealing with
the Segal-Shale-Weil representation.
The Segal-Shale-Weil representation is certain unitary representation of $Mp(\mathbb{V},\omega)$ on the Hilbert space of complex valued square Lebesgue integrable functions on $\mathbb{L}.$ We shall denote this representation by $U,$ i.e.,
$$U: Mp(\mathbb{V},\omega) \to \mathcal{U}({\bf L}^2(\mathbb{L})),$$
where $\mathcal{U}({\bf W})$ is the group of unitary operators on a Hilbert space ${\bf W}.$
Let us only mention that the construction of the Segal-Shale-Weil representation is based on the so called Schr\"{o}dinger representation of the Heisenberg algebra of $(\mathbb{V}=\mathbb{L}\oplus \mathbb{L}',\omega)$ and a use of the Stone-von Neuman theorem. For more details, see Weil \cite{Weil} or Habermann, Habermann \cite{HH}. The Segal-Shale-Weil representation is not irreducible and decomposes into a direct sum of two irreducible modules
${\bf L}^2(\mathbb{L})_+$ and ${\bf L}^2(\mathbb{L})_-$  of even and odd complex square Lebesgue integrable functions on the Lagrangian subspace $\mathbb{L},$ respectively.

For some technical reasons, we shall use the  notion of the so called symplectic spinors. Let us denote the minimal globalization of the Harish-Chandra underlying $(\mathfrak{g},\tilde{K})$-module of the Segal-Shale-Weil representation ${\bf L}^2(\mathbb{L})$ by $\bf S$, i.e., ${\bf S}= mg(HC({\bf L}^2(\mathbb{L}))).$ Let us denote the resulting representation of $\tilde{G}$ on $\bf S$ by $meta$ and call it the {\it metaplectic representation}, i.e.,
$$meta: Mp(\mathbb{V},\omega) \to \hbox{Aut}({\bf S}).$$ We shall call the elements of $\bf S$ {\it symplectic spinors}. According to the decomposition
${\bf L}^2(\mathbb{L})={\bf L}^2(\mathbb{L})_+\oplus{\bf L}^2(\mathbb{L})_-,$ we have also ${\bf S} =  {\bf S}_+\oplus {\bf S}_-.$ Thus $\bf S$ splits into a direct sum of two irreducible representations ${\bf S}_+$ and ${\bf S}_-.$
It is well known that as  $\mathfrak{g}^{\mathbb{C}}$-modules, the modules $s({\bf S}_+)$ and $HC({\bf S}_-)$ are isomorphic to $L(\lambda_0)$ and $L(\lambda_1),$ respectively, where $\lambda_0:=-\frac{1}{2}\varpi_l$ and $\lambda_1=\varpi_{l-1}-\frac{3}{2}\varpi_l.$
   Thus, we have $\mathbb{S}_+\simeq L(\lambda_0)$ and $\mathbb{S}_- \simeq L(\lambda_1).$ For more details, see Kashiwara, Vergne \cite{KV}.

   It is well known that $\mathbb{S}$ is equivalent (as a $\mathfrak{g}^{\mathbb{C}}$-module) to the so called Dixmier realization of the complex symplectic Lie algebra $\mathfrak{g}^{\mathbb{C}}$.   Because of this equivalence, we shall denote this realization by $meta,$ too. See also the section 2.1. for the  remarks on the notation.  
The Dixmier realization $meta$ of $\mathfrak{sp}(\mathbb{V}^{\mathbb{C}},\omega)$   is an injective Lie algebra homomorphism $meta: \mathfrak{sp}(\mathbb{V}^{\mathbb{C}},\omega)\to \hbox{End}_{\mathbb{C}}(\mathbb{C}[x^1,\ldots, x^l])$ of the complex symplectic Lie algebra into the Lie algebra of endomorphisms of polynomials.
It is given as follows (cf. Britten, Hooper, Lemire \cite{BHL}).
\begin{eqnarray}
meta(X_{\epsilon_i-\epsilon_{i+1}})&:=&x_{l-i}\frac{\partial}{\partial x^{l-i+1}}, \quad i=1,\ldots, l-1 \\
meta(X_{-(\epsilon_i-\epsilon_{i+1})})&:=&x_{l-i+1}\frac{\partial}{\partial x^{l-i}}, \quad i=1,\ldots, l-1\\
meta(X_{2\epsilon_l})&:=&-\frac{1}{2}\frac{\partial^2}{\partial (x^1)^2} \\
meta(X_{-2\epsilon_l})&:=&\frac{1}{2}(x^1)^2.
\end{eqnarray}

 Thus, we have $\mathbb{S}\simeq \mathbb{C}[x^1,\ldots, x^l]$ as $\mathfrak{g}^{\mathbb{C}}$-modules.
Moreover, we have that $meta$ restricted to the invariant subspaces of even polynomials is equivalent to $\mathbb{S}_+$ 
and $meta$ restricted to the odd polynomials is equivalent to $\mathbb{S}_-.$
Using this identification, the highest weight vector of $\mathbb{S}_+$ is $1 \in \mathbb{C}[x^1,\ldots, x^l]$ 
 and the highest weight vector of $\mathbb{S}_-$ is $x^1\in \mathbb{C}[x^1,\ldots, x^l].$
For more information see, e.g., Britten, Hooper, Lemire \cite{BHL} and Kashiwara, Vergne \cite{KV}.

In order to derive the studied Howe type duality, we shall also need the symplectic Clifford multiplication $\mathbb{V}\times \bf{S} \to \bf{S}$, which enables us to multiply symplectic spinors by  vectors from $\mathbb{V}.$ It is given by the following prescription
$$(e_i.f)(x)=\frac{\partial f}{\partial x^i}(x),$$
$$(e_{i+l}.f)(x)=\imath x^i f(x),\, i=1,\ldots, l,$$ where $x \in \mathbb{L}, f\in {\bf S}$
and it is extended linearly to the whole $\mathbb{V}.$

This multiplication is basically the Schr\"{o}dinger 
 quantization prescription.
Let us only mention that the differentiation and the multiplication in the definition of the symplectic Clifford multiplication
make sense because of the interpretation of the minimal globalization, see, e.g., Schmid \cite{Schmid}.

{\bf Lemma 4:} For $v,w \in \mathbb{V}$ and $s\in {\bf S},$ the following relation 
$$v.w.s-w.v.s=\imath \omega(v,w)s$$
holds.

{\it Proof.} See Habermann, Habermann \cite{HH}. $\Box$

In the next lemma, certain $\tilde{G}$-equivariance of the symplectic Clifford multiplication is described.

{\bf Lemma 5:} For $g\in Mp(\mathbb{V},\omega),$ $v\in \mathbb{V}$ and $s\in {\bf  S},$ we have $$meta(g)(v.s)=(\lambda(g)v).meta(g)s.$$

{\it Proof.} See Habermann, Habermann \cite{HH}. $\Box$

\subsection{Decomposition of symplectic spinor valued forms}

We shall need the following

{\bf Lemma 6:} Let $\mathbb{V}^{\mathbb{C}}$ be the defining representation of $\mathfrak{g}^{\mathbb{C}},$ then
$$\bigwedge^{r}(\mathbb{V}^*)^{\mathbb{C}} \simeq \bigoplus_{j=0}^{\lfloor r/2 \rfloor}(\omega^{\wedge j} \wedge F(\varpi_{r - 2j}))$$ for $r=0,\ldots, l,$
where $\lfloor q \rfloor$ is the lower integral part of an element $q \in \mathbb{R}$.

{\it Proof.} See, e.g., Corollary 5.1.9. pp. 237 and Theorem 5.1.8.(3) pp. 236 in Goodman and Wallach \cite{GW}. $\Box$

In the previous lemma, the symbol $F(\omega_0)$ denotes the trivial representation of $\mathfrak{sp}(\mathbb{V}^{\mathbb{C}},\omega).$

Let us remark that $\bigwedge^i\mathbb{V}^{\mathbb{C}}\simeq \bigwedge^i(\mathbb{V}^{\mathbb{C}})^*\simeq \bigwedge^{2l-i}(\mathbb{V}^{\mathbb{C}})^*$ as
$\mathfrak{g}^{\mathbb{C}}$-modules because the symplectic form $\omega$ gives a $\mathfrak{g}^{\mathbb{C}}$-isomorphism $\mathbb{V}^{\mathbb{C}}\simeq (\mathbb{V}^{\mathbb{C}})^*$ and because of $Sp(\mathbb{V}
^{\mathbb{C}},\omega)\subseteq SL(\mathbb{V}^{\mathbb{C}}).$

Next, we shall display the decomposition of the tensor product of a finite dimensional $\mathfrak{g}^{\mathbb{C}}$-module and the respective basic symplectic spinor modules $\mathbb{S}_{\pm},$ which was done by  Britten, Hooper and Lemire in \cite{BHL}.
Before doing so, let us introduce the following set. For a dominant integral weight $\lambda=\sum_{j=1}^l\lambda_j \varpi_j$ and $i=0,1,$ we set

\begin{eqnarray*}T^{i}_{\lambda}:=\{\mu \in (\mathfrak{h}^*)^{\mathbb{C}} &|&  \lambda - \mu =:\sum_{j=1}^l d_j\epsilon_j;\, d_j+\delta_{l,j}\delta_{1,i}\in \mathbb{N}_0, j=1,\ldots, l;
 0\leq d_j \leq  \lambda_j,\\
&& j=1,\ldots, l-1; \,
 0 \leq d_l + \delta_{1,i} \leq 2\lambda_l+1 \, \hbox{and} \, \sum_{j=1}^l d_j  \hbox{ is even}\}.
\end{eqnarray*}

{\bf Theorem 7:}
For $i=0,1$ and an integral dominant weight $\lambda$ wr. to $(\mathfrak{h}^{\mathbb{C}},\Phi^+),$ we have
$$F(\lambda) \otimes L(\lambda_i) \simeq  \bigoplus_{\kappa \in T_{\lambda}^i \cap \Pi(\lambda)}L(\lambda_i + \kappa).$$

{\it Proof.} See Britten, Hooper, Lemire \cite{BHL}. $\Box$

Let us remark that there is a misprint in the cited theorem in Britten, Hooper, Lemire \cite{BHL} in the case of the tensor product of a finite dimensional symplectic module $F(\lambda)$ and $L(\lambda_1)=\mathbb{S}_-.$

For our convenience, let us introduce a function $sgn:\{+,-\}\to \{0,1\}$ given by the prescription
$sgn(+):=0$ and $sgn(-):=1.$  The following modules will be used in the decomposition theorem.
$$\mathbb{E}_{ij}^{\pm}:=L(\underbrace{\frac{1}{2},\ldots,\frac{1}{2}}_{j},\underbrace{-\frac{1}{2},\ldots,
-\frac{1}{2}}_{l-j-1},-1+\frac{1}{2}(-1)^{i+j+sgn(\pm)})$$ for $i=0,\ldots,l-1,$   $j=0,\ldots, i$ and $i=l,$   $j=0,\ldots, l-1.$ For $i=j=l,$
we set $\mathbb{E}^+_{ll}:=L(\frac{1}{2}\ldots \frac{1}{2})$ and $\mathbb{E}^-_{ll}:=L(\frac{1}{2}\ldots \frac{1}{2} \frac{5}{2}).$
For $i=l+1,\ldots,2l$ and $j=0,\ldots, 2l-i,$ we set $\mathbb{E}_{ij}^{\pm}:=\mathbb{E}_{(2l-i)j}^{\pm}.$ 

In order to write the results in a short form, let us introduce the following set of pairs of non-negative integers.
$$\Xi:=\{(i,j)| i=0,\ldots,l; j=0,\ldots, i\}\cup \{(i,j)|
i=l+1,\ldots,2l; j=0,\ldots, 2l-i\}\,.
$$ In the canonical coordinate system of $\mathbb{R}^2,$ the set $\Xi$ has a shape of a triangle, see Figure 1. bellow. For $i=0,\ldots, l,$ let us denote
the number $i$ by $m_i$ and for $i=l+1,\ldots, 2l,$ let us denote the number  $2l-i$ by $m_i.$ Thus we can write
$\Xi:=\{(i,j)|i=0,\ldots, 2l; j= 0,\ldots, m_i\}.$ For a convenience, let us set $\mathbb{E}_{ij}^{\pm}:=0$ for $(i,j)\in \mathbb{Z}^2-\Xi.$

In the next lemma, we use two theorems of Goodman, Wallach \cite{GW}, in order to get an information about the set of all weights of the fundamental representations of the complex symplectic Lie algebra $\mathfrak{sp}(\mathbb{V}^{\mathbb{C}},\omega).$ This will enable 
us to use the Theorem 7, where the set $\Pi(\lambda)$ of all weights occurs.

{\bf Lemma 8:} For $r=1,\ldots, l,$
we have $$\Pi(\omega_r) \supseteq \{\sum_{s=1}^r\pm\epsilon^{i_s}|1\leq i_1 < \ldots < i_r  \leq l\}.$$

{\it Proof.} According to the Theorem  5.1.11. pp. 237 and Theorem 6.1.8.(3) pp. 236 in Goodman, Wallach \cite{GW}, the
$\mathfrak{g}^{\mathbb{C}}$-module $F(\omega_r)$ is isomorphic to a linear span
of isotropic $j$-vectors in $\mathbb{V}^{\mathbb{C}}$, i.e., multi vectors $v=f_1\wedge\ldots\wedge f_r,$ where $\omega(f_i,f_j)=0$ for $i,j=1,\ldots, l.$ Second, it is easy to realize that one can choose the Cartan subalgebra $\mathfrak{h}^{\mathbb{C}}$ of $\mathfrak{g}^{\mathbb{C}}$ in a way that the following is true: For $i=1,\ldots, l,$ the basis vector $e_i \in \mathbb{V}^{\mathbb{C}}$ is a weight vector of the weight $\epsilon_i$ and the weight vector $e_{i+l}$ is a  weight vector of the weight $-\epsilon_i,$ both of the defining representation $\mathfrak{g}^{\mathbb{C}}$ on $\mathbb{V}^{\mathbb{C}}.$ Using this, the result follows.
$\Box$

For another realization of $F(\omega_j)$ and for a dimension formula for these representations, see Goodman, Wallach \cite{GW}. Let us also remark that one could prove a more precise result as that one written in the Lemma 8.  Namely, one could prove an equality instead of the inclusion.  But in this case, one should either add the zero element to the set at the right hand side, or not change this set, depending on certain circumstances. We shall only use the information written in the Lemma 8.

Now, we shall introduce a module $\bf W$ we shall be dealing with. As a vector space,
${\bf W}:=\bigwedge^{\bullet}\mathbb{V}^{*}\otimes {\bf S}.$ We assume $\bf W$ to be equipped with the Grothendieck  tensor product topology. For details on this topology, see Treves \cite{Tr} and Vogan \cite{Vogan}. The representation $ \rho: \tilde{G} \to \mbox{Aut}(\bf W)$ of $\tilde{G}$ on $\bf W$ we are interested in, is defined by the following prescription

 $$\rho(g)(\alpha\otimes s):=\lambda^{*\wedge r}(g)\alpha \otimes meta(g)s$$ for
 $g\in \tilde{G},$ $\alpha\in \bigwedge^r \mathbb{V}^*$ and $s\in {\bf S}.$

Now, we can state the decomposition theorem. The proof of this theorem is based on a direct application of the Lemma 6, the Theorem 7 and the Lemma 8.

{\bf Theorem 9:} For $i=0,\ldots,2l,$ the following decomposition into irreducible $\mathfrak{g}^{\mathbb{C}}$-modules
$$\bigwedge^i (\mathbb{V}^{*})^{\mathbb{C}}\otimes \mathbb{S}_{\pm}=\bigoplus_{(i,j)\in\Xi} \mathbb{E}_{ij}^{\pm}$$
holds.

{\it Proof.} Let us remark that because of the fact $\mathbb{V}\simeq \mathbb{V}^*$ as $\mathfrak{g}$-modules, we will write  $\mathbb{V}$  instead of $\mathbb{V}^*$ in this proof.

Using the Lemma 6, we obtain that for $i=2k,$ $k\in \mathbb{N}_0,$ we have
\begin{eqnarray}
\bigwedge^i \mathbb{V}^{\mathbb{C}}\otimes \mathbb{S}_{\pm}=(F(\omega_0)\oplus F(\omega_2)\oplus \ldots \oplus F(\omega_i))\otimes \mathbb{S}_{\pm} \label{1}
\end{eqnarray}
For $i=2k+1,$ $k\in\mathbb{N}_0,$ we have
\begin{eqnarray}
\bigwedge^i \mathbb{V}^{\mathbb{C}}\otimes \mathbb{S}_{\pm}=(F(\omega_1)\oplus F(\omega_3)\oplus \ldots \oplus F(\omega_i))\otimes \mathbb{S}_{\pm} \label{2}
\end{eqnarray}

We shall be considering the mentioned tensor products  for $i=0,\ldots,l$ only, because the result for  $i=l+1,\ldots,2l$ follows from that one for $i=0,\ldots,l$ immediately
due the  $\mathfrak{g}^{\mathbb{C}}$-isomorphism $\bigwedge^i \mathbb{V}^{\mathbb{C}}\otimes \mathbb{S}_{\pm} \simeq \bigwedge^{2l-i} \mathbb{V}^{\mathbb{C}}\otimes \mathbb{S}_{\pm}$   and the definition of $\mathbb{E}_{ij}^{\pm}$ for
$i=l+1,\ldots,2l$ and $j=0,\ldots, 2l-i.$

Let us consider the tensor products by $\mathbb{S}_+$ and $\mathbb{S}_-$ separately.

\begin{itemize}
\item[1.] First, let us consider the tensor product $\bigwedge^i \mathbb{V}^{\mathbb{C}}\otimes \mathbb{S}_+$
for $i=0,\ldots,l.$

Using the Lemma 8 and the Theorem 7, we can easily compute that for $j=1,\ldots, l,$
$T_{\omega_j}^0=\{\epsilon_1+\ldots+\epsilon_j, \epsilon_1 + \ldots +\epsilon_{j-1}- \epsilon_l\} \subseteq \Pi(\omega_j)$ and
$F(\omega_j)\otimes \mathbb{S}_+=L(\underbrace{\frac{1}{2},\ldots,\frac{1}{2}}_j,\underbrace{-\frac{1}{2},\ldots,-\frac{1}{2}}_{l-j})\oplus
L(\underbrace{\frac{1}{2},\ldots,\frac{1}{2}}_{j-1},\underbrace{-\frac{1}{2},\ldots,-\frac{1}{2}}_{l-j},-\frac{3}{2}),$
using the relation $\omega_j=\sum_{i=1}^j\epsilon_i.$
Summing up these terms according to (\ref{1}) and (\ref{2}), we obtain the statement of the theorem for both of the cases $i$  is odd and $i$
is even.

\item[2.]  Now, let us consider the tensor product $\bigwedge^i \mathbb{V}^{\mathbb{C}}\otimes \mathbb{S}_-$ for $i=0,\ldots,l.$

 Using the Lemma 8, we can easily compute that for $j=1,\ldots, l-1,$ we have
$T_{\omega_j}^1=\{\epsilon_1+\ldots+\epsilon_j,\epsilon_1+\ldots+\epsilon_{j-1}+\epsilon_l\}\subseteq \Pi(\omega_j)$ and $T^1_{\omega_l}=\{\epsilon_1+\ldots + \epsilon_l, \epsilon_1+\ldots + \epsilon_{l-1}-\epsilon_l\}\subseteq \Pi(\omega_l).$ 
Therefore using the Theorem 7, we have 
$F(\omega_j)\otimes \mathbb{S}_-=L(\underbrace{\frac{1}{2},\ldots,\frac{1}{2}}_{j-1},\underbrace{-\frac{1}{2},\ldots,-\frac{1}{2}}_{l-j+1})\oplus L(\underbrace{\frac{1}{2},\ldots,\frac{1}{2}}_{j},\underbrace{-\frac{1}{2},\ldots,-\frac{1}{2}}_{l-j-1},-\frac{3}{2})$ for $j=1,\ldots, l-1.$ For $j=l,$ we obtain
$F(\omega_l)\otimes \mathbb{S}_+=L(\frac{1}{2} \ldots \frac{1}{2} -\frac{3}{2})\oplus L(\frac{1}{2} \ldots \frac{1}{2}-\frac{5}{2})$ again using the Theorem 7.
Summing up these terms according to (\ref{1}) and (\ref{2}), we obtain the statement of the theorem for both cases
$i$ is odd and $i$ is even.
\end{itemize}
$\Box$

In the next figure (Figure 1.), one can see the decomposition structure of $\bigwedge^{\bullet}(\mathbb{V}^{*})^{\mathbb{C}}\otimes \mathbb{S}_+$ in the case of $l=3.$ For $i=0,\ldots, 6,$ the $i^{th}$ column constitutes of  modules in which the $\mathbb{S}_+$-valued exterior forms of form degree $i$ decompose.

$$\xymatrix{
\mathbb{E}^+_{00}  &\mathbb{E}_{10}^+  &\mathbb{E}_{20}^+  &\mathbb{E}_{30}^+  &\mathbb{E}_{40}^+  &\mathbb{E}_{50}^+  &\mathbb{E}_{60}^+ \\
&\mathbb{E}_{11}^+   &\mathbb{E}_{21}^+ &  \mathbb{E}_{31}^+   & \mathbb{E}_{41}^+ &\mathbb{E}_{51}^+  & \\
& & \mathbb{E}_{22}^+   & \mathbb{E}_{32}^+   & \mathbb{E}_{42}^+ &&\\
&&& \mathbb{E}_{33}^+ &&&}
$$

\centerline{Figure 1.}

\section{Representation of $\mathfrak{osp}(1|2)$ on $\bigwedge^{\bullet}\mathbb{V}^* \otimes {\bf S}$}
In this section, we shall introduce certain simple super Lie algebra and its representation in order to describe  the commutant algebra $\mbox{End}_{\tilde{G}}(\bigwedge^{\bullet}\mathbb{V}^* \otimes {\bf S})$  in the next section.

Let us say a few words about the notation concerning super Lie algebras we shall be using. For each super Lie algebra $\mathfrak{f}$, let us denote the even and the odd part by $\mathfrak{f}_0$ and $\mathfrak{f}_1,$ respectively. Sometimes, the generators of $\mathfrak{f}_0$ and $\mathfrak{f}_1$ are
called bosonic and fermionic generators, respectively.
The super Lie bracket of two homogeneous elements $u, v$ will be denoted by $[u,v]$ if and only if at least one of them is from the even part $\mathfrak{f}_0.$ In the other cases, we shall denote it by $\{u,v\}.$

Let $\mathfrak{g}'$ be the simple super Lie algebra $\mathfrak{osp}(1|2).$
Using the usual definition of $\mathfrak{osp}(1|2),$ there exists a basis $\{ h, e^+, e^-, f^+, f^-\}$ of $\mathfrak{g}',$ such that $\{e^+,h,e^-\}$ spans the even part $\mathfrak{g}'_0,$ $\{f^+,f^-\}$ spans the odd part $\mathfrak{g}'_1$ and
the only nonzero relations among the basis elements are

\begin{equation}
[h,e^{\pm}]= \pm e^{\pm} \qquad [e^+,e^-]=2h \label{he}
\end{equation}

\begin{equation}
[h,f^{\pm}]=\pm \frac{1}{2}f^{\pm} \qquad \{f^{+},f^-\}=\frac{1}{2}h \label{hf}
\end{equation}

\begin{equation}
[e^{\pm},f^{\mp}]=-f^{\pm} \qquad \{f^{\pm},f^{\pm}\}=\pm \frac{1}{2}e^{\pm} \label{ef}
\end{equation}
 
Now, we shall introduce the following five operators which will be then used to define a representation of $\mathfrak{g}'=\mathfrak{osp}(1|2).$

For   $r=0,\ldots, 2l$ and $\alpha \otimes s \in \bigwedge^r
\mathbb{V}^*\otimes {\bf S},$  we set
$$F^+:\bigwedge ^{r}\mathbb{V}^*\otimes {\bf S} \to
\bigwedge^{r+1}\mathbb{V}^*\otimes {\bf S},\, F^+(\alpha \otimes
s):=\frac{\imath}{2}\sum_{i=1}^{2l}\epsilon^i\wedge \alpha \otimes e_i.s,$$
$$F^-:\bigwedge ^{r}\mathbb{V}^*\otimes {\bf S} \to \bigwedge ^{r-1}\mathbb{V}^*\otimes {\bf S}, \,
F^-(\alpha \otimes s):=\frac{1}{2}\sum_{i=1}^{2l}\iota_{\check{e}_i}\alpha \otimes
e_i.s$$
and extend them linearly.

Next, we shall define the operators $H, E^+$ and $E^-$ in the following way.
For $r=0,\ldots, 2l,$ we set
$$H:\bigwedge ^{r}\mathbb{V}^*\otimes {\bf S} \to \bigwedge ^{r}\mathbb{V}^*\otimes {\bf S},\, H=\{F^+,F^-\}$$
$$E^{\pm}:\bigwedge ^{r}\mathbb{V}^*\otimes {\bf S} \to \bigwedge ^{r\pm 2}\mathbb{V}^*\otimes {\bf S}, \,
 E^{\pm}:=\pm 2 \{F^{\pm},F^{\pm}\}$$ where $\{,\}$ means the anticommutator in the associative algebra $\mbox{End}({\bf W}).$

Now, we introduce a $\mathbb{Z}_2$-graded structure on the vector space ${\bf W}.$
Let us set ${\bf W}_0 :=(\bigoplus_{i = 0}^{l}\bigwedge^{2i}\mathbb{V}^*) \otimes {\bf S}_+$ and ${\bf W}_1 :=(\bigoplus_{i=0}^{l-1}\bigwedge^{2i+1}\mathbb{V}^*) \otimes {\bf S}.$ Now, we can say that $\mbox{End}(\bigwedge^{\bullet}\mathbb{V}^{*} \otimes {\bf S})$ has the inherited super Lie algebra structure from the super vector space structure on ${\bf W}={\bf W}_0\oplus {\bf W}_1.$ This means that as a super vector  space, we have $\hbox{End}({\bf W})=\hbox{End}_0({\bf W})\oplus \hbox{End}_1({\bf W}),$ where
$\hbox{End}_0({\bf W})=\hbox{End}({\bf W}_0)\oplus \hbox{End}({\bf W}_1)$ and $\hbox{End}_1({\bf W})=\hbox{Hom}({\bf W}_0,{\bf W}_1)\oplus \hbox{Hom}({\bf W}_1,{\bf W}_0)$ as vector spaces. (The introduced grading on $\hbox{End}({\bf W})$ together with composition of endomorphisms make it a super associative algebra. The super Lie algebra structure is the inherited structure from the super associative structure.) 

Let us consider the following
mapping $$\sigma: \mathfrak{osp}(1|2) \to  \hbox{End}(\bigwedge^{\bullet}(\mathbb{V}^*)^{\mathbb{C}}\otimes {\bf S})$$ defined by $$\sigma(e^{\pm}):=E^{\pm},$$
$$\sigma(f^{\pm}):=F^{\pm}$$ and $$\sigma(h):=H$$
and extend it linearly to the whole $\mathfrak{g}'=\mathfrak{osp}(1|2).$

In what follows, we shall prove that $\sigma$ is a super Lie algebra representation (Theorem 11).  Before doing so, let us make some preliminary computations in order to have simpler formulas for the
introduced mappings $E^{\pm}, H$ and $F^{\pm}.$

Using the antisymmetry of $\wedge,$ commutator relation of the symplectic Clifford multiplication (Lemma 4), renumbering of indices, and the antisymmetry of $\omega_{ij},$ we get
\begin{eqnarray*}
E^{+}(\alpha\otimes s)  &=& 2 \{F^+,F^+\}(\alpha \otimes s)\\
   &=& -\sum_{i,j=1}^{2l}\epsilon^i\wedge \epsilon^j\wedge \alpha \otimes e_i.e_j.s \\
   &=& \sum_{i,j=1}^{2l}\epsilon^j\wedge\epsilon^i\wedge \alpha \otimes (e_j.e_i. - \imath \omega_{ij})s\\
   &=& \sum_{i,j=1}^{2l}\epsilon^i\wedge\epsilon^j \wedge \alpha \otimes (e_i.e_j. + \imath \omega_{ij})s \\
   &=&-E^{+}(\alpha \otimes s) +  \imath\sum_{i,j=1}^{2l}\omega_{ij}\epsilon^i\wedge \epsilon^j\wedge \alpha \otimes s,
\end{eqnarray*}
where $\alpha \otimes s \in \bigwedge^r(\mathbb{V}^*)^{\mathbb{C}}\otimes {\bf S}$ and $i=0,\ldots,2l.$

Comparing the left hand and right hand side, we get
\begin{eqnarray}
E^+(\alpha \otimes s)=\frac{\imath}{2}\sum_{i,j=1}^{2l}  \omega_{ij}\epsilon^i\wedge\epsilon^j \wedge \alpha \otimes s.
\end{eqnarray}

Similarly, one  computes
\begin{eqnarray}
E^{-}(\alpha \otimes s)=  \frac{\imath}{2}\sum_{i,j=1}^{2l} \omega_{ij}\iota_{\check{e}_i}\iota_{\check{e}_j}\alpha \otimes s,
\end{eqnarray} where
$\alpha \otimes s \in \bigwedge^r(\mathbb{V}^*)^{\mathbb{C}}\otimes {\bf S}$ and $i=0,\ldots,2l.$

A similar computation for the operator $H$ is a bit more complicated, and it is treated in the next

{\bf Lemma 10:} Let $(\mathbb{V},\omega)$ be a $2l$ dimensional
symplectic vector space. Then for   $r=0,\ldots, 2l,$ we have
$$H_{|\bigwedge^r(\mathbb{V}^*)^{\mathbb{C}} \otimes {\bf S}}=\frac{1}{2}(r-l)\hbox{Id}_{|\bigwedge^r(\mathbb{V}^*)^{\mathbb{C}}\otimes {\bf S}}.$$

{\it Proof.} We shall compute the value of $H$ in three items.
\begin{itemize}
\item[(i)] First, let us prove that for $r=0,\ldots, 2l $ and $\alpha \in
\bigwedge^r(\mathbb{V}^*)^{\mathbb{C}},$ we have
\begin{eqnarray}
\sum_{i,j=1}^{2l}\omega_{ij}\epsilon^i\wedge
\iota_{\check{e}_j}\alpha = r \alpha.
\end{eqnarray}
It is sufficient to prove the statement for homogeneous basis
elements. We proceed by induction.

I. For $r=1,$ we have $\sum_{i,j=1}^{2l}\omega_{ij}\epsilon^i
\wedge\iota_{\check{e}_j}\epsilon^k=
\sum_{i,j=1}^{2l}\omega_{ij}\epsilon^i
\omega^{jk}=\delta_i^k\epsilon^i=\epsilon^k$ for each $k=1,\ldots,
2l.$

II. For $r=0,\ldots, 2l,$ $\alpha \in \bigwedge^r(\mathbb{V}^*)^{\mathbb{C}}$ and
$k=1,\ldots, 2l,$ we have
$\sum_{i,j=1}^{2l}\omega_{ij}\epsilon^i\wedge\iota_{\check{e}_j}(\epsilon^k\wedge\alpha)=
\sum_{i,j=1}^{2l}\omega_{ij}\epsilon^i\wedge\iota_{\check{e}_j}\epsilon^k\wedge\alpha
- \sum_{i,j=1}^{2l}\omega_{ij}\epsilon^i\wedge \epsilon^k \wedge
\iota_{\check{e}_j}\alpha = \sum_{i,j=1}^{2l} \omega_{ij}
\omega^{jk}\\ \epsilon^i \wedge \alpha + \sum_{i,j=1}^{2l}
\omega_{ij} \epsilon^k \wedge \epsilon^i \wedge
\iota_{\check{e}_j}\alpha= \epsilon^k\wedge \alpha +
r\epsilon^k\wedge \alpha = (r+1)\epsilon^k\wedge \alpha,$ where we
have used the induction hypothesis in the second last equation.

\item[(ii)] Now, let us prove that for each  $s\in {\bf S},$
\begin{eqnarray}
A:=\sum_{i,j=1}^{2l}\omega^{ij} e_i.e_j.s =  \imath l s.
\end{eqnarray}
 Using the Lemma 4 about the commutator of the symplectic Clifford
multiplication, we can write
$A=\sum_{i,j=1}^{2l}\omega^{ij}e_i.e_j.s=\sum_{i,j=1}^{2l}\omega^{ij}(e_j.e_i.-\imath\omega_{ij})s=
\sum_{i,j=1}^{2l}\omega^{ji}e_i.e_j.s+\sum_{i=1}^{2l}\imath\delta^{i}_is=
-\sum_{i,j=1}^{2l} \omega^{ij}e_i.e_j.s+ 2\imath l s= -A + 2\imath l
s.$ Comparing the left hand and right hand side, we get $A=\imath l
s.$
\item[(iii)] Now, we shall compute the value of  $H.$
For $r=0,\ldots, 2l$ and $\alpha \otimes s \in
\bigwedge^r(\mathbb{V}^*)^{\mathbb{C}}\otimes {\bf S},$ we have
$H(\alpha\otimes s)=2(F^+F^-(\alpha\otimes s)+F^-F^+(\alpha
\otimes s))=
\frac{\imath}{2}(\sum_{i,j=1}^{2l}\epsilon^i\wedge\iota_{\check{e}_j}\alpha\otimes
e_i.e_j.s + \sum_{i,j=1}^{2l}\iota_{\check{e}_i}(\epsilon^j\wedge
\alpha)\otimes e_i.e_j.s)=
\frac{\imath}{2}(\sum_{i,j=1}^{2l}\epsilon^i\wedge\iota_{\check{e}_j}\alpha\otimes
e_i.e_j.s + \sum_{i,j=1}^{2l}\iota_{\check{e}_i}\epsilon^j\alpha
\otimes e_i.e_j.s- \sum_{i,j=1}^{2l}\epsilon^j \wedge
\iota_{\check{e}_i}\alpha\otimes e_i.e_j.s)=
\frac{\imath}{2}(\sum_{i,j=1}^{2l}\epsilon^j\wedge\iota_{\check{e}_i}\alpha\otimes
e_j.e_i.s + \sum_{i,j=1}^{2l}\alpha \otimes \omega^{ij}e_i.e_j.s-
\sum_{i,j=1}^{2l}\epsilon^j \wedge (\iota_{\check{e}_i}
\alpha)\otimes (e_j.e_i.s - \imath\omega_{ij}s))=
\frac{\imath}{2}(\sum_{i,j=1}^{2l}\alpha \otimes
\omega^{ij}e_i.e_j.s -\sum_{i,j=1}^{2l} \imath \omega_{ji}
\epsilon^j \wedge \iota_{\check{e}_i}\alpha\otimes s)=
\frac{\imath}{2}(\imath l\alpha\otimes  s - \imath r \alpha \otimes
s)=\frac{1}{2}(r-l)\alpha \otimes s,$ where we have used the equations
(17) and (18) in the second  last step.
\end{itemize}
$\Box$ 

{\bf Theorem 11:} The mapping
$$\sigma: \mathfrak{osp}(1|2)\to \hbox{End}(\bigwedge^{\bullet}(\mathbb{V}^*)^{\mathbb{C}}\otimes {\bf S})$$
is a super Lie algebra representation.

{\it Proof.} First, it is easy to see that $\sigma(\mathfrak{g}'_i) \subseteq \mbox{End}_i({\bf W}),$ $i=0,1.$

Second, we shall verify that the operators $E^{\pm}, H, F^{\pm}$ satisfy the same relations as that ones for $e^{\pm}, h, f^{\pm}$ written in the rows (\ref{he}), (\ref{hf}), (\ref{ef}).
It is evident from the definitions of $H$ and
$E^{\pm}$ that the second relations of (\ref{hf}) and (14) are satisfied.   Let us only check the $+$ version of the first equation of the row labeled by (\ref{hf}) and the second equation of the row labeled by (\ref{he}). The other relations could be treated in a similar way.

For $r=0,\ldots, 2l$ and $\alpha \otimes s \in \bigwedge^r (\mathbb{V}^*)^{\mathbb{C}} \otimes {\bf S},$ we have
\begin{eqnarray*}
[H,F^{+}](\alpha \otimes s)&=& H F^+ (\alpha \otimes s) - F^+ H(\alpha \otimes s)\\
&=& H (\frac{\imath}{2}\sum_{i=1}^{2l}\epsilon^i\wedge \alpha \otimes e_i .s) - F^+ \frac{1}{2}(r-l)(\alpha \otimes s)\\
&=& \sum_{i=1}^{2l}\left[\frac{1}{2}\frac{\imath}{2} (r+1-l)\epsilon^i \wedge \alpha \otimes e_i.s  - \frac{\imath}{2}\frac{1}{2}(r-l)\epsilon^i \wedge \alpha \otimes e_i.s\right]\\
&=&\frac{\imath}{4} \sum_{i=1}^{2l}\epsilon^i \wedge \alpha \otimes e_i.s = \frac{1}{2}F^+ (\alpha \otimes s).
\end{eqnarray*}
Thus we got the relation in the form of the $+$ version of the first equation written in the row $(\ref{hf})$ as required.

Now, we shall check whether $[E^+,E^-] = 2H.$ In this part, we shall use the Einstein summation convention.
Using the equations (15) and (16), we get for $r= 0,\ldots, 2l$ and $\alpha \otimes s \in \bigwedge^r (\mathbb{V}^*)^{\mathbb{C}}\otimes {\bf S},$
\begin{eqnarray*}
[E^+,E^-](\alpha \otimes s)
&=&E^+(\frac{\imath}{2}\omega_{ij}\iota_{\check{e_i}} \iota_{\check{e_j}} \alpha \otimes s) -
E^-(\frac{\imath}{2}\omega_{ij}\epsilon^i\wedge\epsilon^j\wedge\alpha \otimes s)\\
&=&-\frac{1}{4}\omega_{ij}\omega_{kl}\epsilon^k\wedge \epsilon^l \wedge \iota_{\check{e_i}}\iota_{\check{e_j}} \alpha \otimes s
 +\frac{1}{4}\omega_{kl}\omega_{ij}\iota_{\check{e_k}}\iota_{\check{e_l}}(\epsilon^i\wedge \epsilon^j\wedge \alpha) \otimes s\\
&=& -\frac{1}{4}\omega_{ij}\omega_{kl}\epsilon^k\wedge \epsilon^l \wedge \iota_{\check{e_i}}\iota_{\check{e_j}} \alpha \otimes s  +\frac{1}{4}\omega_{ij}\omega_{kl}\omega^{li}\omega^{kj} \alpha \otimes s\\
   &&- \frac{1}{4}\omega_{ij}\omega_{kl}\omega^{li}\epsilon^j\wedge  \iota_{\check{e_k}}\alpha\otimes s
 -\frac{1}{4}\omega_{ij}\omega_{kl}\omega^{lj}\omega^{ki}\alpha \otimes s\\
&&  +\frac{1}{4}\omega_{ij}\omega_{kl}\omega^{lj} \epsilon^i\wedge \iota_{\check{e_k}}\alpha \otimes s
 + \frac{1}{4}\omega_{ij}\omega_{kl}\omega^{ki} \epsilon^j \wedge \iota_{\check{e_l}}\alpha \otimes s\\
&&  - \frac{1}{4}\omega_{ij}\omega_{kl}\omega^{kj}\epsilon^i \wedge \iota_{\check{e_l}} \alpha \otimes s
+ \frac{1}{4}\omega_{ij}\omega_{kl}\epsilon^i \wedge \epsilon^j \wedge \iota_{\check{e_k}} \iota_{\check{e_l}} \alpha \otimes s\\
&=&\frac{1}{4}(\omega_{kj}\omega^{kj}\alpha - \omega_{kj}\epsilon^j \wedge \iota_{\check{e_k}}\alpha  -
 \omega_{ij}\omega^{ji} \alpha \\
&&+ \omega_{ik}\epsilon^i \wedge \iota_{\check{e_k}}\alpha -\omega_{lj}\epsilon^j\wedge \iota_{\check{e_l}}\alpha+\omega_{il}\epsilon^i \wedge \iota_{\check{e_l}} \alpha ) \otimes s\\
&=&\frac{1}{4}(-4l\alpha  + 4 \omega_{ij}\epsilon^i \wedge \iota_{\check{e_j}}\alpha)\otimes s\\
&=& (-l + r)\alpha \otimes s = 2H (\alpha \otimes s),
\end{eqnarray*}
where we have used the relation (17) in the second last step.
The other relations can be obtained in a similar way.
$\Box$

{\bf Lemma 12:} The linear mappings $E^{\pm},F^{\pm}$ and $H$ are $\tilde{G}$-equivariant
with respect to the representation $\rho$ of $\tilde{G}$ on ${\bf W}.$

{\it Proof.} Let us prove the equivariance for $F^+, F^-$ and $H$ separately.
\begin{itemize}
\item[(i)] For $g\in \tilde{G},$ $r=0,\ldots, 2l$ and $\alpha\otimes s \in \bigwedge^{r}(\mathbb{V}^*)^{\mathbb{C}}\otimes {\bf S},$ we have
\begin{eqnarray*}
-2\imath \rho(g)  F^{+} (\alpha \otimes s)&=&\rho(g)[\sum_{i=1}^{2l}(\epsilon^{i}\wedge\alpha \otimes e_i.s)]\\
                          &=&\sum_{i=1}^{2l}\lambda(g)^{*}\epsilon^{i}\wedge \lambda(g)^{*}\alpha \otimes meta(g)(e_i.s).
\end{eqnarray*}
Now, we use the equivariance property of the symplectic Clifford multiplication, see Lemma 5.

\begin{eqnarray*}
                          &&\sum_{i=1}^{2l}\lambda(g)^{*}\epsilon^{i}\wedge \lambda(g)^{*}\alpha \otimes meta(g)(e_i.s)\\
                          &=&\sum_{i,j=1}^{2l}\epsilon^j[\lambda(g)^{-1}]_{j}^{\ i}   \wedge \lambda(g)^{* }\alpha \otimes (\lambda (g)e_i). meta(g)\\
                          &=&\sum_{i,j,k=1}^{2l}\epsilon^{j}[\lambda (g)^{-1}]_{j}^{\ i}  \wedge \lambda (g)^{*  } \alpha \otimes [\lambda(g)]^{\ k}_{i}
                          e_k.(meta(g)s)\\
                          &=&\sum_{i,j,k=1}^{2l} [\lambda(g)^{-1}]_{j}^{\ i} [\lambda (g)]^{\ k}_{i}\epsilon^{j}\wedge \lambda(g)^{*  } \alpha \otimes e_k. meta(g)s\\
                          &=&\sum_{j,k=1}^{2l} \epsilon^{j} \delta^k_j \wedge \lambda(g)^{* } \alpha \otimes e_k. meta(g)s
                           = -2\imath F^{+} \rho(g) (\alpha \otimes s).
\end{eqnarray*}
\item[(ii)]
Now, we do a similar computation for $F^{-}.$ For $g \in \tilde{G},$
$r=0,\ldots, 2l$ and  $\alpha\otimes s \in
\bigwedge^{r}(\mathbb{V}^*)^{\mathbb{C}}\otimes {\bf S},$ we obtain
\begin{eqnarray*}
2\rho(g) F^{-} (\alpha \otimes s)&=& 2\rho(g)[\sum_{i=1}^{2l}\iota_{\check{e}_i}\alpha \otimes e_i.s]\\
&=&\sum_{i=1}^{2l} \lambda(g)^{*}\iota_{\check{e}_i}\alpha\otimes meta(g)(e_i.s).\\
\end{eqnarray*}

It is easy to see that
$\lambda(g)^*(\iota_v\alpha)=\iota_{\lambda(g)v}(\lambda(g)^*\alpha)$
for $v\in (\mathbb{V}^*)^{\mathbb{C}}.$ Using this relation, we get
\begin{eqnarray*}
&& \sum_{i=1}^{2l} \lambda(g)^{*}(\iota_{\check{e}_i}\alpha)\otimes meta(g)(e_i.s)\\
&=& \sum_{i,j=1}^{2l}\iota_{\lambda(g)\check{e}_i}(\lambda(g)^{* }\alpha)\otimes [\lambda(g)]^{\ j}_{i}e_j.meta(g)s.\\
\end{eqnarray*}
Using the formula $\check{e}_i=\sum_{k=1}^{2l}\omega^{ik}e_k,$  we can rewrite the preceding equation as
\begin{eqnarray*}
&&\sum_{i,j=1}^{2l} \iota_{\lambda(g)\sum_{k=1}^{2l}\omega^{ik}e_k}(\lambda(g)^{* }\alpha)\otimes [\lambda(g)]^{\ j}_{i}e_j.meta(g)s\\
&=&\sum_{i,j,k,n=1}^{2l}[\lambda(g)]^{\ j}_{i}\omega^{ik}\iota_{[\lambda(g)]^{\ n}_k e_n}(\lambda(g)^{* }\alpha)\otimes e_j.meta(g)s\\
&=&\sum_{i,j,k,n=1}^{2l}[\lambda(g)]^{\ j}_{i}\omega^{ik}[\lambda(g)]^{\ n}_k\iota_{ e_n}(\lambda(g)^{* }\alpha)\otimes e_j.meta(g)s.\\
\end{eqnarray*}
Because $\lambda(g)\in Sp(\mathbb{V},\omega),$ we have $\sum_{i,k=1}^{2l}
[\lambda(g)]^{\ j}_{i}\omega^{ik}[\lambda(g)]^{\ n}_k =\\
\sum_{i,k=1}^{2l}[\lambda(g)^{\perp}]^{j}_{\
i}\omega^{ik}[\lambda(g)]^{\ n}_{k} = \omega^{jn}$ for $j,n =
1,\ldots, 2l.$ Substituting this relation into the previous
computation, we get
\begin{eqnarray*}
&&\sum_{i,j,k, n=1}^{2l}[\lambda(g)]^{\ j}_{i}\omega^{ik}[\lambda(g)]^{\ n}_k\iota_{ e_n}(\lambda(g)^{* }\alpha)\otimes e_j.meta(g)s\\
&=&\sum_{j,n=1}^{2l}\omega^{jn}\iota_{e_n}(\lambda(g)^{* }\alpha)\otimes e_j.meta(g)s\\
&=&\sum_{j=1}^{2l}\iota_{\sum_{n=1}^{2l}\omega^{jn}e_n}(\lambda(g)^{* }\alpha)\otimes e_j.meta(g)s\\
&=&\sum_{j=1}^{2l}\iota_{\check{e}_j}(\lambda(g)^{*}\alpha)\otimes e_j.meta(g)s\\
&=&2F^{-}\rho(g)(\alpha\otimes s).
\end{eqnarray*}

\item[(iii)] $H$ is $\tilde{G}$-equivariant, because $H=2\{F^{+},F^{-}\}$ and $F^{+}$ and $F^{-}$ are $\tilde{G}$-equivariant. The mapping $E^{\pm}
$ is also $\tilde{G}$-equivariant because \\
$E^{\pm}=\pm 2 \{ F^{\pm},F^{\pm}\}$  and $F^{+}$ and $F^{-}$ are $\tilde{G}$-equivariant.
\end{itemize}
$\Box$

Summing up, we have the following

{\bf Corollary 13:} The super Lie algebra representation $$\sigma: \mathfrak{osp}(1|2) \to \mbox{End}(\bigwedge^{\bullet}{(\mathbb{V}}^*)^{\mathbb{C}}\otimes {\bf S})$$ maps the super Lie algebra
$\mathfrak{osp}(1|2)$ into the algebra  $\mbox{End}_{\tilde{G}}(\bigwedge^{\bullet}(\mathbb{V}^{*})^{\mathbb{C}}\otimes {\bf S})$ of $\tilde{G}$-invariants.

{\it Proof.} The fact that $\sigma$ is a super Lie algebra representation 
was proved  in the Theorem 11 and the fact that the image of $\sigma$ lies in the commutant algebra of $\tilde{G}$-invariants 
 was proved
 in the Lemma 12.$\Box$

\section{Computation of $\hbox{End}_{\tilde{G}}(\bigwedge^{\bullet}\mathbb{V}^*\otimes {\bf S})$}

As a short hand, let us set $p_+:=1 \in \mathbb{C}[x^1,\ldots, x^l]$ and $p_-:=x^1 \in \mathbb{C}[x^1,\ldots, x^l].$ 
For $(i,j) \in \Xi$ and $i-j$ even, let us define a vector $$v_{ij}^{\pm}:=\omega^{\wedge \lfloor \frac{i-j}{2}\rfloor}\wedge \epsilon^{l+1}\wedge \ldots \wedge \epsilon^{l+j} \otimes p_{\pm}.$$
In the case $i-j$ is odd, let us take 
$$v_{ij}^{+}:=\omega^{\wedge \lfloor \frac{i-j}{2} \rfloor}\wedge \epsilon^l \wedge \epsilon^{l+1}\wedge \ldots\wedge \epsilon^{l+j}\otimes p_{+} \mbox{ and}$$
$$v_{ij}^-:=\omega^{\wedge \lfloor \frac{i-j}{2}\rfloor}\wedge \epsilon^{l+1}\wedge\ldots \wedge \epsilon^{l+j}\wedge \epsilon^{2l}\otimes p_-.$$

We will use this vectors in order to prove that restrictions of $F^{\pm}$ to the modules $\mathbb{E}_{ij}^{\pm}$ for suitable $i,j$ are non-zero. These restrictions will be then used to prove that the image of the representation $\sigma$ and two additional operators generates the algebra $\hbox{End}_{\tilde{G}}(\bigwedge^{\bullet}\mathbb{V}^*\otimes {\bf S}).$ 

{\bf Lemma 14:} For $(i,j)\in \Xi,$ the vector
$v_{ij}^{\pm}$ is a highest weight vector of $\mathbb{E}_{ij}^{\pm}.$

{\it Proof.} Using the the formulas (1), (2), (6), (8) for computing the action of $\mathfrak{g}^{\mathbb{C}}$ on the vector $v^{\pm}_{ij},$ we see that  $v_{ij}^{\pm}$ is a maximal vector in $\bigwedge^i(\mathbb{V}^*)^{\mathbb{C}}\otimes \mathbb{S}_{\pm},$ $(i,j)\in \Xi.$ Because
$\bigwedge^i(\mathbb{V}^*)^{\mathbb{C}}\otimes \mathbb{S}_{\pm}$ decomposes into a finite set of mutually non-equivalent irreducible representations (Theorem 9), there is only one (up to a nonzero complex multiple) maximal vector of a given weight. Using the formulas (4), (5) and the fact that the weight of $p_{\pm}$ is $\lambda_{sgn(\pm)},$  one can compute the weight of  $v_{ij}^{\pm}$ for $(i,j)\in \Xi.$ One finds out that the weight of $v_{ij}^{\pm}$ is precisely the highest weight of $\mathbb{E}^{\pm}_{ij}.$ (See the definition of $\mathbb{E}_{ij}^{\pm}.)$ $\Box$


For our convenience, let us define the following two subsets $\Xi_+$ and $\Xi_-$ of the set $\Xi.$ We set $\Xi_+:=\Xi-\{(j,2l-j)| j=l,\ldots, 2l\}$ and $\Xi_-:=\Xi-\{(j,j)|
 j=0,\ldots, l\}.$
The sets $\Xi_+$ and $\Xi_-$ could be represented by the left hand side and the right hand side of the "triangle" $\Xi,$ respectively. See the Figure 1.
Let us consider an element $\beta \in \bigwedge^{\bullet}(\mathbb{V}^*)^{\mathbb{C}}$ and write it in the form
$\beta=\sum_{I\subseteq \{1,\ldots, 2l\}}\beta_I\epsilon^I,$ where $\epsilon^I=\epsilon^{i_1}\wedge \ldots \wedge \epsilon^{i_k}$ if $I=(i_1,\ldots, i_k),$ $\beta_I\in \mathbb{C}$ and $1\leq k \leq 2l.$ We say, that the expression for $\beta$ is in a standard form if we are summing
over the indices $I$ of  the form $I=(i_1,\ldots, i_k)$ for $1\leq i_1 < \ldots < i_k\leq 2l,$ $1\leq k \leq 2l.$

{\bf Lemma 15:} For each $(i,j)\in \Xi,$ we have
$$F^{+}_{|\mathbb{E}_{ij}^{\pm}}: \mathbb{E}_{ij}^{\pm}\stackrel{\sim}{\rightarrow}
\mathbb{E}_{i+1,j}^{\mp} \qquad \mbox{ if } (i,j) \in \Xi_+,$$
$$F^{-}_{|\mathbb{E}_{ij}^{\pm}}: \mathbb{E}_{ij}^{\pm}\stackrel{\sim}{\rightarrow}
\mathbb{E}_{i-1,j}^{\mp} \qquad \mbox{ if } (i,j)\in \Xi_-.$$

 {\it Proof.} We do this proof for the operator $F^{+}$ only. The case of the
operator $F^-$ can be treated in a similar way.
\begin{itemize}
\item[1.)]First, we prove the statement about the target space. From
the definition of $F^+,$ we have that for $ (i,j) \in \Xi,$ the image $\mbox{Im}(F^+_{|\mathbb{E}_{ij}^{\pm}})\subseteq \bigwedge^{i+1}(\mathbb{V}^*)^{\mathbb{C}} \otimes \mathbb{S}_{\mp}$ as a $\mathfrak{g}^{\mathbb{C}}$-submodule. 
Due to the structure of
$\bigwedge^{i+1}(\mathbb{V}^*)^{\mathbb{C}} \otimes \mathbb{S}_{\mp},$ we see that
$\mbox{Im}(F^+_{|\mathbb{E}_{ij}^{\pm}}) \subseteq \mathbb{E}_{i+1,j}^{\mp}$
because $F^{+}$ is  $\mathfrak{g}^{\mathbb{C}}$-equivariant (Lemma 12) and
$\mathbb{E}_{i+1,j}^{\mp}$ is the only irreducible submodule
of  $\bigwedge^{i+1}(\mathbb{V}^*)^{\mathbb{C}}\otimes \mathbb{S}_{\mp}$ isomorphic to
$\mathbb{E}_{ij}^{\pm}$ (we suppose $(i,j)\in \Xi_+$).
\item[2.)] Let us consider the case $\bigwedge^{\bullet}(\mathbb{V}^*)^{\mathbb{C}}\otimes \mathbb{S}_+$ only and take a pair of non-negative integers $(i,j)\in \Xi_+.$
 For a fixed $i-j$ even, let us compute 
\begin{eqnarray}
F^+(v_{ij}^+)&=&\frac{\imath}{2}\omega^{\wedge\lfloor \frac{i-j}{2} \rfloor}\wedge\sum_{r=1}^{2l}\epsilon^r\wedge \epsilon^{l+1}\wedge \ldots \wedge\epsilon^{l+j}\otimes e_r.1\nonumber \\
&=& -\frac{1}{2}\omega^{\wedge\lfloor \frac{i-j}{2} \rfloor}\wedge\sum_{r=l+1}^{2l} \epsilon^r \wedge \epsilon^{l+1}\wedge \ldots \wedge\epsilon^{l+j}\otimes x^r\nonumber \\
&=&-\frac{1}{2}\omega^{\wedge\lfloor \frac{i-j}{2} \rfloor}\wedge\sum_{r=l+j+1}^{2l}\epsilon^r \wedge \epsilon^{l+1}\wedge \ldots \wedge\epsilon^{l+j} \otimes x^r
\end{eqnarray}
Now, we shall prove that the last written expression is non-zero provided $(i,j) \in \Xi_+.$
Let us write the term $\omega^{\wedge \lfloor \frac{i-j}{2}\rfloor}$ in the standard form. All summands in $\omega^{\wedge\lfloor \frac{i-j}{2} \rfloor}$ containing $\epsilon^t$ for some 
$t\in \{l+1,\ldots, l+j\}$ contribute to zero and could be neglected. Let us denote the resulting set of the considered summands by $M.$ 
To each member of this set, we can associate an unordered $2\lfloor \frac{i-j}{2}\rfloor$-tuple consisting of the labels of the expressions 
which are exteriorly multiplied in the chosen summand. Let us denote the set of such tuples by $\mathcal{M}.$ For example, if $\epsilon^1\wedge \epsilon^{2l}\in M,$ then $(1,2l)\in \mathcal{M}.$ Now, it is sufficient to prove that there is an $r \in \{l+j+1, \ldots, 2l \}$ and an element $s \in \mathcal{M}$ such that $r \notin s.$ Indeed, if this is true, we know that there is a non-zero summand in $F^+(v_{ij}^+)$ and this summand does not cancel with another summand which is of the form $\alpha\otimes x^r$ for some  exterior form $\alpha.$ (If it cancels, then the element $\omega^{\wedge \lfloor \frac{1-j}{2}\rfloor}$ would not be in the standard form). Let us prove the mentioned existence of $r \in \{l+j+1,\ldots, 2l\}.$ If $i=j$ we are done, see (19). Suppose $i>j.$
From the structure of $\omega$ and the construction of $\mathcal{M},$
we know that there is an element $s \in \mathcal{M}$ of the form $s=(j+1,\ldots,j+\lfloor \frac{i-j}{2} \rfloor, j+l+1,\ldots, j+\lfloor \frac{i-j}{2} \rfloor + l).$ If we prove that the largest element of $s$ is smaller then $r=2l,$ we are done.
This leads to the inequality $j+\lfloor \frac{i-j}{2} \rfloor <l.$ 
Because $i-j$ is even, the last written inequality translates into $i+j<2l$ which is true because of $(i,j)\in \Xi_+.$ Using the fact that $\mathbb{E}_{ij}^{\pm}$ are irreducible and $F^{\pm}$ are $\tilde{G}$-equivariant (Lemma 12), the result follows.
The results for the complementary cases ($i-j$ odd, $F^-$ and $\bigwedge^{\bullet}(\mathbb{V}^*)^{\mathbb{C}} \otimes \mathbb{S}_-$) 
could be treated in a similar way.$\Box$
\end{itemize}

{\bf Remark:} It is  easy to see   that $F^-$ is zero when restricted to 
$\mathbb{E}_{ii}^{\pm},$ $i=0,\ldots,l.$  Namely, we know that $F^-$ is lowering the form degree by one, it is a $\mathfrak{g}^{\mathbb{C}}$-equivariant map and  there is no submodule of the module of symplectic spinor valued extrerior forms of form degree  $i-1$ isomorphic to $\mathbb{E}_{ii}^{\pm}.$ 
A similar discussion could be done for $F^+.$

Let us define the following operators $R_{\pm}: \bigwedge^{\bullet}\mathbb{V}^* \otimes \mathbb{S}\to \bigwedge^{\bullet}\mathbb{V}^* \otimes \mathbb{S}_{\pm}$ by the formula
$R_{\pm}(\alpha \otimes s):=\alpha \otimes s_{\pm},$ where $s=(s_{+},s_-)$ according to the decomposition $\mathbb{ S} =  \mathbb{S}_+\oplus \mathbb{S}_-,$ and extend it linearly.
For $(i,j)\in \Xi,$ let us denote the projection operators from the space $\bigwedge^{\bullet}(\mathbb{V}^*)^{\mathbb{C}}\otimes \mathbb{S}_{\pm}$ to the submodule
$\mathbb{E}_{ij}^{\pm}$ by $S_{ij}^{\pm},$ i.e.,
$$S_{ij}^{\pm}:\bigwedge^{\bullet}(\mathbb{V}^*)^{\mathbb{C}} \otimes \mathbb{S}_{\pm}\to \mathbb{E}_{ij}^{\pm} \subseteq \bigwedge^i (\mathbb{V}^*)^{\mathbb{C}} \otimes \mathbb{S}_{\pm}.$$ 
(By a projection operator $S$ onto a vector space, we mean a surjective operator satisfying $S^2=S,$ i.e., we do not demand any hermicity condition in particular.)
Due to the structure of the decomposition of $\bigwedge^i(\mathbb{V}^*)^{\mathbb{C}}\otimes \mathbb{S}_{\pm}$ (see Theorem 9), the operators $S^{\pm}_{ij}$ are uniquely defined.

{\bf Lemma 16:} For each $(i,j)\in \Xi,$ the projection operators
$S_{ij}^{\pm}$ are elements of the associative algebra generated by
$F^{\pm}$ and $R_{\pm}$ as an associative algebra over $\mathbb{C}.$

{\it Proof.} For each $i = 0, \ldots, 2l,$ let us  define the following operators $S_i^{\pm}:
\bigwedge^{\bullet}(\mathbb{V}^*)^{\mathbb{C}} \otimes \mathbb{S}_{\pm} \to
\bigwedge^{i}(\mathbb{V}^*)^{\mathbb{C}}\otimes \mathbb{S}_{\pm}$ by the formula
$$S_i^{\pm}:=\left( \prod_{r=0, r\neq i}^{2l}\frac{2H-r+l}{i-r}\right) R_{\pm}.$$
Using the Lemma 10, we see that the image of each $S_i^{\pm}$ lies in the
prescribed space. Recall that $H$ could be expressed by $F^{+}$ and
$F^-$ and thus for $i=0,\ldots, 2l,$ the projections $S_i^{\pm}$ lie in
the algebra described in the formulation of the lemma. Now let us fix an element $i\in \{0,\ldots, 2l\}.$ We prove that for
each $j$ such that $(i,j)\in \Xi,$ the projections
$S_{ij}^{\pm}:\bigwedge^{i}(\mathbb{V}^*)^{\mathbb{C}}\otimes \mathbb{S}_{\pm} \to \mathbb{E}_{ij}^{\pm}$ 
could be written as linear combinations of
compositions of $F^{\pm}$ and $R_{\pm}.$ We proceed by
induction.
\begin{itemize}
\item[I.] For $j=0,$ we can define ${S_{i0}^{\pm}}':=(F^+)^{i}(F^-)^{i}.$ Using
the fact that applying the $F^-$ (or $F^{-}$) lowers (or rises) the form
degree by $1,$ we see that ${S_{i0}^{\pm}}': \bigwedge^{i}(\mathbb{V}^*)^{\mathbb{C}}\otimes
\mathbb{S}_{\pm}\to \mathbb{E}_{i0}^{\pm}$ (see the Figure 1.). Using the Schur lemma (Theorem
3), we see that there exists a complex number $\lambda_{i0}\in
\mathbb{C}$ such that ${S_{i0}^{\pm}}''_{|\mathbb{E}^{\pm}_{i0}}=\lambda_{i0}\mbox{Id}_{|\mathbb{E}_{i0}^{\pm}}.$ 
Due to the Lemma 15, we know that
$\lambda_{i0}\neq 0.$ Thus defining
$S_{i0}^{\pm}:=\frac{1}{\lambda_{i0}}{S_{i0}^{\pm}}'\circ S_i^{\pm},$ we get the projection onto
$\mathbb{E}_{i0}^{\pm}$ expressed as linear combinations of compositions of
$F^{\pm}$ and $R_{\pm}.$

\item[II.]Let us suppose that for $k=0,\ldots, j,$ the operators $S_{ik}^{\pm}$
could be written as linear combinations of compositions of the
operators $F^{\pm}$ and $R_{\pm}.$ Now, we shall use the
operators $S_{i0}^{\pm},\ldots, S_{ij}^{\pm}$ in order to define $S_{i,j+1}^{\pm}.$
Let us take $\xi \in \bigwedge^{i}(\mathbb{V}^*)^{\mathbb{C}}\otimes \mathbb{S}_{\pm}$
and form an element $\zeta:= {S_{i,j+1}^{\pm}}'\xi:=\xi -
\sum_{k=0}^jS_{ik}^{\pm}\xi \in \bigoplus_{k=j+1}^{m_i} \mathbb{E}_{ik}^{\pm}$ (see the Figure 1.). Now, form an element
$\zeta':={S_{i,j+1}^{\pm}}''\zeta:=(F^+)^{i-j}(F^-)^{i-j}\zeta.$ In the
same way as in the item I., we conclude that $\zeta' \in \mathbb{E}_{i,j+1}^{\pm}.$ 
 Using the Schur lemma (Theorem 3) in the case of $S_{i,j+1 |\mathbb{E}_{i,j+1}^{\pm}}'': \mathbb{E}_{i,j+1}^{\pm} \to \mathbb{E}_{i,j+1}^{\pm},$ we
conclude that there is a complex number $\lambda_{i,j+1} \in
\mathbb{C}$ such that ${S_{i,j+1|\mathbb{E}_{ij}^{\pm}}^{\pm}}''=\lambda_{i,j+1}
\mbox{Id}_{|\mathbb{E}_{i,j+1}^{\pm}}.$  Due to the Lemma 15, we know that
$\lambda_{i,j+1}\neq 0.$ Now define
$S_{i,j+1}^{\pm}:=\frac{1}{\lambda_{i,j+1}}{S_{i,j+1}^{\pm}}''\circ S_i^{\pm}.$ One easily sees, that
$S_{i,j+1}^{\pm}$ is the desired projection. $\Box$
\end{itemize}

Let us denote the complex associative algebra generated over $\mathbb{C}$ by $F^+,$ $F^-$ and the projection operators $R_{\pm}$ by $C.$ 
Now, we prove that the operators $\sigma(\mathfrak{g}')$ and the operators
$R_{\pm}$ give enough $\tilde{G}$-invariants in the sense of the following

{\bf Theorem 17:} We have the following associative algebra  isomorphism 
between the algebra of $\tilde{G}$-invariants $\hbox{End}_{\tilde{G}}(\bigwedge^{\bullet}(\mathbb{V}^*)^{\mathbb{C}} 
\otimes {{\bf S}})$  and the algebra $C,$ i.e.,
$$\mbox{End}_{\tilde{G}}(\bigwedge^{\bullet}(\mathbb{V}^*)^{\mathbb{C}}\otimes {\bf S}) \simeq C.$$

{\it Proof.} The proof is based on the decomposition of
$\bigwedge^{\bullet}{(\mathbb{V}^{*})}^{\mathbb{C}}\otimes \mathbb{S}_{\pm}$ into
irreducible summands written in the Theorem 9. 
Due to the Corollary 12, we know that $C \subseteq
\hbox{End}_{\tilde{G}}(\bigwedge^{\bullet}(\mathbb{V}^*)^{\mathbb{C}}\otimes {\bf S}).$ Thus, we should prove the opposite inclusion
$\hbox{End}_{\tilde{G}}(\bigwedge^{\bullet} (\mathbb{V}^*)^{\mathbb{C}}\otimes {\bf S})
\subseteq C.$ Clearly, we have $\mbox{End}_{\mathfrak{g}^{\mathbb{C}}}(\mathbb{W}^{\mathbb{C}})=
\mbox{End}_{\mathfrak{g}^{\mathbb{C}}}(\bigwedge^{\bullet}{(\mathbb{V}^*)}^{\mathbb{C}}\otimes \mathbb{S}_+)\oplus
\mbox{End}_{\mathfrak{g}^{\mathbb{C}}}(\bigwedge^{\bullet}{(\mathbb{V}^*)}^{\mathbb{C}}\otimes \mathbb{S}_-)\oplus
\mbox{Hom}_{\mathfrak{g}^{\mathbb{C}}}(\bigwedge^{\bullet}{(\mathbb{V}^*)}^{\mathbb{C}}\otimes \mathbb{S}_+,\bigwedge^{\bullet}{(\mathbb{V}^*)}^{\mathbb{C}}\otimes \mathbb{S}_-) \oplus \\
\oplus \mbox{Hom}_{\mathfrak{g}^{\mathbb{C}}}(\bigwedge^{\bullet}{(\mathbb{V}^*)}^{\mathbb{C}}\otimes \mathbb{S}_-,\bigwedge^{\bullet}{(\mathbb{V}^*)}^{\mathbb{C}}\otimes \mathbb{S}_+).$ Let us only
prove $\mbox{Hom}_{\mathfrak{g}^{\mathbb{C}}}(\bigwedge^{\bullet}{(\mathbb{V}^*)}^{\mathbb{C}}\otimes \mathbb{S}_+,\bigwedge^{\bullet}{(\mathbb{V}^*)}^{\mathbb{C}}\otimes \mathbb{S}_-)\subseteq C.$
The remaining three cases could be handled in a similar way.
Due to the structure of the decomposition of $\bigwedge^{\bullet}(\mathbb{V}^{*})^{\mathbb{C}}\otimes \mathbb{S}$  described in the Theorem 9, we see that
$\hbox{Hom}_{\mathfrak{g}^{\mathbb{C}}}(\bigwedge^{\bullet}{(\mathbb{V}^*)}^{\mathbb{C}}\otimes \mathbb{S}_+,
\bigwedge^{\bullet}{(\mathbb{V}^*)}^{\mathbb{C}}\otimes \mathbb{S}_-)\simeq \bigoplus_{(i,j), (k,m)\in \Xi}
\hbox{Hom}_{\mathfrak{g}^{\mathbb{C}}}(\mathbb{E}_{ij}^+, \mathbb{E}_{km}^-).$   
Regarding the highest weights of the irreducible $\mathfrak{g}^{\mathbb{C}}$-modules $\mathbb{E}_{ij}^{+}$ and $\mathbb{E}_{km}^{-}$ and using the 
 the remark bellow the Schur lemma (Theorem 3), we see that the following is true. Namely, if there exists a non trivial $\mathfrak{g}^{\mathbb{C}}$-equivariant map between $\mathbb{E}_{ij}^+$ and $\mathbb{E}_{km}^{-},$ then $j=m$ and $i\equiv k+1 \, \hbox{mod} \, 2.$
Due to the Lemma 15 and the remark bellow this lemma,
we know that $(F^{+})^{1+2r}S_{ij}^{+}: \mathbb{E}_{ij}^+
\to \mathbb{E}_{i+2r+1,j}^-$  for each $r\in \mathbb{N}_0$ and $(i,j)\in \Xi.$
Due to the Lemma 15 again, we know that the mapping is nontrivial whenever $(i,j)\in \Xi_+$ and $(i+2r+1,j) \in \Xi.$
Thus $\mathbb{C}(F^{+})^{1+2r}S_{ij}^+=\mbox{Hom}_{\mathfrak{g}^{\mathbb{C}}}(\mathbb{E}_{i,j}^+,\mathbb{E}_{i+1+2r,j}^-)$ due to the remark bellow the Schur lemma (Theorem 3).
A similar reasoning could be done for $F^-$ (case $r\in \{-1,-2,\ldots\}$).
Using the fact that $S_{ij}^{+}$ could be written only by using  linear combinations of compositions
 of the operators $F^{\pm}$ and $R_{\pm}$ (Lemma 16), we have
 $C\simeq  \mbox{End}_{\mathfrak{g}^{\mathbb{C}}}(\bigwedge^{\bullet}(\mathbb{V}^*)^{\mathbb{C}}\otimes \mathbb{S}).$ 
 Because the considered representation of $\mathfrak{g}^{\mathbb{C}}$ is complex and $\mathfrak{g}$ is a real form of $\mathfrak{g}^{\mathbb{C}},$
 we see that $C \simeq \mbox{End}_{\mathfrak{g}}(\bigwedge^{\bullet}(\mathbb{V}^*)^{\mathbb{C}}\otimes \mathbb{S}).$ (We suppose $\bigwedge^{\bullet}(\mathbb{V}^*)^{\mathbb{C}}\otimes \mathbb{S}$ to be a complex representation.)
 Because $\tilde{K}$ is connected,
  we have $C \simeq \mbox{End}_{(\mathfrak{g},\tilde{K})}(\bigwedge^{\bullet}(\mathbb{V}^*)^{\mathbb{C}}\otimes S)$ due to the Theorem 2.
   Using the Theorem 1 about the minimal globalization functor, we have the isomorphism
 $C \simeq \mbox{End}_{\tilde{G}}(\bigwedge^{\bullet}(\mathbb{V}^*)^{\mathbb{C}}\otimes {\bf S}).$
$\Box$


{\bf Remark:} Let us remark that more conceptually, we could have defined the projections $S_{ij}^{\pm}$ onto the submodules $\mathbb{E}_{ij}^{\pm}$ using the Casimir operator $Cas \in \mathcal{U}(\mathfrak{osp}(1|2))$ of the super Lie algebra $\mathfrak{osp}(1|2).$ (The symbol $\mathcal{U}(\mathfrak{f})$ is used for the universal enveloping algebra of the super Lie algebra $\mathfrak{f}.$)



{\bf Remark:} In the case of a spin group acting on (orthogonal) spinor valued exterior forms, the Howe
dual could be described by a use of the Lie algebra $\mathfrak{sl}(2,\mathbb{C}).$ See the Introduction for a reference. The algebra
$\mathfrak{osp}(1|2)$ is often considered as the super symmetric
analogue of the algebra $\mathfrak{sl}(2,\mathbb{C}).$ Thus, we see
that there is a "symmetry" in changing the symmetric form
(orthogonal case) into an antisymmetric one (symplectic case).

Summing up, we can say a bit informally (see the Introduction) that the ortho-symplectic  super Lie algebra $\mathfrak{osp}(1|2)$ is the Howe dual partner of the metaplectic group in the studied case.

The use of the Howe dual in representation theory is not purposeless. It is often used in order to decompose a module into irreducible pieces and to "hide" their possible  multiplicities. For this purpose, we shall use the super Lie algebra $\mathfrak{osp}(1|2)$ in the next section.

 \section{Two-folded Howe type duality}
 We start defining certain family of finite dimensional  irreducible representations of the super Lie ortho-symplectic algebra $\mathfrak{osp}(1|2).$
   For $j=0,\ldots, l,$ let $\mathbb{G}^j$ be a fixed complex vector space of complex dimension $2l-2j+1$ and consider a basis $\{f_{i}\}_{i=j}^{2l-j}$ of $\mathbb{G}^j.$ (The numbering of the basis elements in the above form is set for convenience.)
The super vector space structure on $\mathbb{G}^j$ is defined as follows. For $j=0,\ldots, l,$ the even part $(\mathbb{G}^j)_0$ is spanned by the basis vectors with an even index, i.e., $(\mathbb{G}^j)_0:=\hbox{Span}_{\mathbb{C}}(\{f_i| i \in \{j,\ldots, 2l-j\}\cap 2\mathbb{N}_0  \}).$ Complementary, we define $(\mathbb{G}^j)_1:=\hbox{Span}_{\mathbb{C}}(\{f_i| i \in \{j,\ldots, 2l-j\} \cap (2\mathbb{N}_0 + 1)\}).$
For our convenience, we set $f_{r}:=0$ for $r \in \mathbb{Z}-\{j,\ldots, 2l-j\}.$ 
We shall not  denote the dependence of this basis on the number $j$ explicitly. We hope this will cause no confusion. As a short hand for each $(i,j)\in \Xi,$ we introduce the rational numbers
$$A(l,i,j):=\frac{(-1)^{i-j}+1}{16}(i-j) + \frac{(-1)^{i-j+1}+1}{16}(i + j - 2l - 1).$$

For $j=0,\ldots, l,$ we define
the mentioned representation 
$$\sigma_j: \mathfrak{osp}(1|2)\to \hbox{End}(\mathbb{G}^j)$$
\begin{eqnarray*}
\sigma_j(f^+)(f_i)&:=&A(l,i+1,j)f_{i+1},\, i= j,\ldots, 2l-j,\\
\sigma_j(f^-)(f_i)&:=&f_{i-1}, \, i=j,\ldots, 2l-j,\\
\sigma_j(h)&:=&2\{\sigma_j(f^+), \sigma_j(f^{-})\}, \\
\sigma_j(e^{\pm})&:=&\pm 2\{\sigma_j(f^{\pm}),\sigma_j(f^{\pm})\}.
\end{eqnarray*}
  
We have the following  

 {\bf Lemma 18:} For $j=0,\ldots, l,$ the mapping $\sigma_j: \mathfrak{osp}(1|2)\to \hbox{End}(\mathbb{G}^j)$ is an irreducible representation of the super Lie algebra $\mathfrak{osp}(1|2).$ 

{\it Proof.}
First we prove that for $j=0,\ldots, l,$ the mapping  $\sigma_j$ is a super Lie algebra representation of the algebra $\mathfrak{osp}(1|2).$
It is easy to see that  whereas the even part of $\mathfrak{g}'$ acts by transforming the even part of $\mathbb{G}^j$ into itself and the odd part into itself as well,  the odd part of $\mathfrak{g}'$ acts by interchanging the mentioned two parts of $\mathbb{G}^j.$ Now, we should check weather the relations (12), (13) and (14) are preserved by the mapping $\sigma_j$ for $j=0,\ldots, l.$
The last relation in the row (13) is satisfied by the definition of $\sigma_j(h),$ and the last and the second last relation in the row (14) are  satisfied by definition (of $\sigma_j(e^+)$ and $\sigma_j(f^-)$) as well.
To prove the other relations is straightforward. We shall only check that the last relation in the row (12) holds.  
We shall consider the following two cases.
\begin{itemize}
\item[1.)] The case $i-j$ is even. Chose $f_i$ for $i=j,\ldots, 2l-j.$
The left hand side of the last relation in the row $(12)$ reads.
\begin{eqnarray*}
&&[\sigma_j(e^+)\sigma_j(e^-)-\sigma_j(e^-)\sigma_j(e^+)]f_i\\
&&=16(-\sigma_j(f^+)\sigma_j(f^+)\sigma_j(f^-)\sigma_j(f^-) +
\sigma_j(f^-)\sigma_j(f^-)\sigma_j(f^+)\sigma_j(f^+))f_i\\
&&=16(-\sigma_j(f^+)\sigma_j(f^+)f_{i-2}+\sigma_j(f^-)\sigma_j(f^-)A(l,i+2,j)A(l,i+1,j)f_{i+2})\\
&& = 16(-A(l,i,j)A(l,i-1,j)f_i+A(l,i+2,j)A(l,i+1,j)f_i)\\
&& = 16(\frac{2.2}{16.16})[-(i-j)(i+j-2l-2)\\
&& + (i+2-j)(i+1+j-2l-1)]f_i\\
&& = (i - l)f_i.
\end{eqnarray*}
(The reader could check that for $i=j,j+1,2l-j-1, 2l-j,$ the equality between the third and the fourth row of the preceding computation is also valid.
In the mentioned cases, the appropriate coefficients $A(l,r,s)$ are zero.)
Now, let us compute the right hand side.
\begin{eqnarray*}
&&   2\sigma_j(h)f_i=4(\sigma_j(f^+)\sigma_j(f^-)+\sigma_j(f^-)\sigma_j(f^+))f_i\\
&& = 4(\sigma_j(f^+)f_{i-1}+A(l,i+1,j)f_{i})\\
&& = 4[\frac{2}{16}(i-j)f_i+\frac{2}{16}(i+1+j-2l-1)f_i]\\
&& = (i - l)f_{i}.
\end{eqnarray*}
(In a similar way as before, we can treat the cases $i=j,j+1, 2l-j-1, 2l-j.)$ 
Thus the left hand side of $[\sigma_{j}(e^+),\sigma_j(e^-)]=\sigma_j(h)$ equals the right hand side of this equation.

\item[2.)] The case $i-j$ is odd could be handled in a similar way and will be omitted.
\end{itemize}

Now, we prove that for $j=0,\ldots, l,$ the representation $\sigma_j$ is irreducible.
Let us suppose, there exists a non-trivial proper invariant subspace $X \subseteq \mathbb{G}^j$ of dimension $0<k<2l-2j+1.$
Let us chose a basis $\{v_j,\ldots, v_{j+k-1}\}$ of $X.$  Let $Y:=\hbox{Span}_{\mathbb{C}}(\{f_j,\ldots, f_{j+k-1}\})$ and define a mapping $T: X \to Y$ by the prescription
$Y(v_i):=f_i,$ $i=j,\ldots, j+k-1$ and extend it first linearly and then to an automorphism $\tilde{T}$ of $\mathbb{G}^j.$
Now, for each $j=0,\ldots,l$, consider a representation $\tilde{\sigma}_j: \mathfrak{osp}(1|2) \to \hbox{End}(\mathbb{G}^j)$ given by $\tilde{\sigma}_j(X):=\tilde{T}^{-1}\sigma_j(X)\tilde{T}$ for each $X\in \mathfrak{osp}(1|2).$ This representation is clearly equivalent to $\sigma_j.$
Applying $\sigma_j(f^+)$ to $f_{j+k-1}\in Y,$ we get a vector $0 \neq f \notin Y.$ We have $\tilde{\sigma_j}(f^+)v_{j+k-1}=\tilde{T}^{-1}\sigma_j(f^+)\tilde{T} v_{j+k -1}
=\tilde{T}^{-1}\sigma_j(f^+)f_{j+k-1} = \tilde{T}^{-1}f \notin X,$  i.e.,  $X$ is not invariant.
Thus we see, that $\tilde{\sigma}_j$ is irreducible for $j =0 ,\ldots, l.$ Because of the equivalence
of $\tilde{\sigma}_j$ and $\sigma_j,$ the lemma follows.
$\Box$

{\bf Lemma 19:} For each $k \in \mathbb{N}_0$ and $i=0,\ldots, 2l,$ we have
\begin{eqnarray*}
&&(F^-)^k F^+=\\
&&(-1)^k F^+(F^-)^k + \left[\frac{(-1)^k+1}{16}k+\frac{(-1)^{k+1}+1}{16}(2i-2l-k+1)\right](F^{-})^{k-1},
\end{eqnarray*}
when acting on $\bigwedge^i(\mathbb{V}^*)^{\mathbb{C}}\otimes {\bf S}.$

{\it Proof.} We will not write explicitly that we are considering the action of the considered operators on the space $\bigwedge^{i}(\mathbb{V}^*)^{\mathbb{C}}\otimes {\bf S}$ and proceed by induction.
\begin{itemize}
\item[I.] For $k=0$ the lemma holds obviously.
\item[II.] 
\begin{itemize}
\item[a.]We suppose that the lemma holds for an even integer $k\in \mathbb{N}_0.$
We have 
\begin{eqnarray*}
(F^-)^{k+1}F^+&=& F^-(F^-)^kF^+\\
              &=& F^-[F^+(F^-)^k+\frac{k}{16}((-1)^k+1)(F^-)^{k-1}]\\
              &=& -(F^+)(F^-)^{k+1}+\frac{1}{2}H(F^-)^k+2\frac{k}{16}(F^{-})^k\\
              &=& -F^+(F^-)^{k+1}+\frac{1}{4}(i-k-l)(F^{-})^{k}+\frac{k}{8}(F^{-})^k\\
              & & -F^+(F^-)^{k+1}+\frac{2}{16}(2i-(k+1)+1-2l)(F^{-})^k,
\end{eqnarray*}
where we have used the induction hypothesis, commutation relation $2\{F^{+}, F^{-}\}=H$
and the Lemma 10 on the value of $H.$ The last written expression coincides with that one of statement of the lemma for $k+1$ odd.
\item[b.] Now, suppose $k$ is odd.
We have
\begin{eqnarray*}
&&(F^-)^{k+1}F^{+}=F^-(F^-)^{k}F^{+}\\
&=&F^-[-F^+(F^-)^{k}+\frac{(-1)^{k+1}+1}{16}(2i-2l-k+1)(F^{-})^k]\\
&=&+F^+(F^-)^{k+1}-\frac{1}{2}H(F^{-})^k+\frac{1}{8}(2i-2l-k+1)(F^{-})^k\\  
 &=&F^+(F^-)^{k+1}-\frac{1}{8}(2i-2k-2l)(F^-)^k+\\
&&+\frac{1}{8}(2i-2l-k+1)(F^{-})^k =F^+(F^-)^{k+1}-\frac{2}{16}(k+1)(F^{-})^k,
\end{eqnarray*}               
where again we have used the same tools as in the previous item.
\end{itemize}
\end{itemize}
$\Box$

For $(i,j)\in \Xi,$ the vector space ${\bf E}_{ij}^{\pm}$ is defined by the following standard procedure.
The highest weight $\mathfrak{g}^{\mathbb{C}}$-modules $\mathbb{E}_{ij}^{\pm}$ were already defined. The $\tilde{K}$-module structure is inherited from the $\tilde{K}$-module structure on $\bigwedge^{\bullet}(\mathbb{V}^*)^{\mathbb{C}}\otimes {\bf S}.$ Finally ${\bf E}_{ij}^{\pm}:=mg(E_{ij}^{\pm}),$ where $E_{ij}^{\pm}$ is mentioned  the $(\mathfrak{g},\tilde{K})$-module.

Let us introduce the following mapping $Sgn: \{\pm\} \times \mathbb{N}_0 \to \{\pm\}$ given by the prescription $Sgn(\pm, 2k):=\pm$ and $Sgn(\pm,2k+1)=\mp$ for each $k\in \mathbb{N}_0.$
Now for $(i,j) \in \Xi,$ let us define a mapping
$\psi_{ij}^{\pm}: {\bf E}_{ij}^{\pm} \to {\bf E}_{jj}^{Sgn(\pm,i-j)} \otimes \mathbb{G}^j$
by the formula 
$$\psi_{ij}^{\pm} v: = (F^{-})^{i-j}v\otimes f_i,$$ for
an element $v \in {\bf E}_{ij}^{\pm}.$

Let us define a family $\{\rho_j^{\pm}|j=0,\ldots, l\}$ of representation of the group $\tilde{G}=Mp(\mathbb{V},\omega)$ on ${\bf E}_{jj}^{\pm}$ simply by a restriction of the representation $\rho$, i.e.,  
$$\rho_j^{\pm}: \tilde{G} \to \hbox{Aut}({\bf E}_{jj}^{\pm})$$
$\rho_j(g)v:=\rho(g)v$ for $j=0,\ldots, l$ $g\in \tilde{G}$ and $v \in {\bf E}_{jj}^{\pm}.$

In the next theorem, the (2-folded) Howe type duality is stated.

{\bf Theorem 20:} We have the following decomposition
$$\bigwedge^{\bullet}(\mathbb{V}^*)^{\mathbb{C}}\otimes {\bf S} \simeq \bigoplus_{j=0}^{l}
({\bf E}_{jj}^+ \oplus {\bf E}_{jj}^-)\otimes \mathbb{G}^j$$ as an
$(Mp(\mathbb{V},\omega) \times \mathfrak{osp}(1|2))$-module.

{\it Proof.} 
Let us recall that at the left hand side, we are considering the representation $L:=\rho \hat{\otimes} \sigma$ of $\tilde{G}\times \mathfrak{g}'$ and at the right hand side, we are considering the representation $R:=\bigoplus_{j=0}^l(\rho_j^+ \oplus \rho_j^-)\otimes \sigma_j$ of $\tilde{G} \times \mathfrak{g}'.$
Due to the Theorem 9, there exists a $\mathfrak{g}^{\mathbb{C}}$-equivariant isomorphism
$\psi: \bigwedge^{\bullet}(\mathbb{V}^*)^{\mathbb{C}}\otimes \mathbb{S} \to \bigoplus_{(i,j)\in \Xi}(\mathbb{E}_{ij}^+\oplus \mathbb{E}_{ij}^-)=:\mathbb{E}.$ The mapping $\psi$ is also  $\mathfrak{g}'$-equivariant due to the definition of the action $\mathfrak{g}'$ on the space $\mathbb{E}$ as we now recall. Actually, we have defined the representation of $\mathfrak{g}'$ on $\mathbb{E}=\bigoplus_{(i,j)\in \Xi}(\mathbb{E}_{ij}^+\oplus \mathbb{E}_{ij}^{-})$ by transporting the $\mathfrak{g}'$-module structure of $\bigwedge^{\bullet}(\mathbb{V}^*)^{\mathbb{C}}\otimes \mathbb{S}$ by the isomorphism $\psi$ to get the $\mathfrak{g}'$-module structure on  $\mathbb{E}.$ Thus we have to prove that for each $(i,j)\in \Xi,$ the mapping $\psi_{ij}^{\pm}: \mathbb{E}_{ij}^{\pm} \to \mathbb{E}_{jj}^{Sgn(\pm,i-j)}\otimes \mathbb{G}^j$ is $(\mathfrak{g}^{\mathbb{C}}\times \mathfrak{g}')$-equivariant. The $\mathfrak{g}^{\mathbb{C}}$-equivariance follows easily because $F^-$ in the definition of $\psi_{ij}^{\pm}$ commutes with the representation $\rho$ of $\tilde{G}$, see the Lemma 12. We should check the $\mathfrak{g}'$-equivariance.
For each $(i,j)\in \Xi$ and $v \in \mathbb{E}_{ij}^{\pm},$ consider 
$\psi_{ij}^{\pm}F^- v=(F^-)^{i-1-j}F^-v \otimes f_{i-1}=(F^-)^{i-j} v \otimes f_{i-1}.$
On the other hand, we have
$F^-\psi_{ij}v=F^-((F^{-})^{i-j}v\otimes f_i)=(F^{-})^{i-j}v\otimes f_{i-1}.$ 
Now, we check the $\mathfrak{g}'$-equivariance in the case of $F^+.$
We shall use the Lemma 19 to compute $\psi_{ij}^{\pm}F^+v=(F^-)^{i+1-j}F^+v\otimes f_{i+1}=[(-1)^{i+1-j}F^+(F^-)^{i+1-j}v+A(l,i+1,j)(F^-)^{i-j}v]\otimes f_{i+1}=
A(l,i+1,j)(F^{-})^{i-j}v \otimes f_{i+1},$ where we have used the fact that $(F^-)^{i+1-j}v=0$ because $v\in \mathbb{E}_{ij}^{\pm}.$
Let us compute $F^+\psi_{ij}^{\pm}v=F^+(F^-)^{i-j}v\otimes A(l,i+1,j)f_{i+1}.$ Thus the equivariance with respect to $F^+$ is proved.
Because the operators $H,$ $E^+$ and $E^-$ are linear combinations of compositions of the operators $F^+$ and $F^-,$ the $\mathfrak{g}'$-equivariance of $\psi_{ij}^{\pm}$ is proved.
Because $\mathfrak{g}$ is a real form of $\mathfrak{g}^{\mathbb{C}}$ and  the considered representation is complex, we have that $\bigwedge^{\bullet}(\mathbb{V}^*)^{\mathbb{C}}\otimes {\bf S} \simeq \bigoplus_{j=0}^{l}
(\mathbb{E}_{jj}^+ \oplus \mathbb{E}_{jj}^-)\otimes \mathbb{G}^j$ as a $(\mathfrak{g}\times \mathfrak{g}')$-module, where the modules $\mathbb{E}_{jj}^+$ and $\mathbb{E}_{jj}^-$ are also irreducible when considered as complex representations of the real Lie algebra $\mathfrak{g}$ for $j=0,\ldots, l.$
Using the fact that $\tilde{K}$ is connected, we get 
$\bigwedge^{\bullet}(\mathbb{V}^*)^{\mathbb{C}}\otimes S \simeq \bigoplus_{j=0}^{l}
[(E_{jj}^+\otimes \mathbb{G}^j) \oplus (E_{jj}^-\otimes \mathbb{G}^j)]$
as a $((\mathfrak{g},\tilde{K})\times \mathfrak{g}')$-module (cf. Theorem 2).
Using the minimal globalization functor and the fact that the Grothendieck tensor product topology is compatible with this functor (see Vogan \cite{Vogan}), we get the globalized situation described in the statement of this lemma.
$\Box$
  
\section{Application: Symplectic Dirac operators and their higher spin analogues.}

In 1995, K. Habermann introduced the so called symplectic Dirac operator using the symplectic spinor bundle introduced by Bertram Kostant in \cite{Kostant}.
The symplectic Dirac operator was introduced mainly in order to derive topological and geometric properties of symplectic manifolds using analytic properties of the introduced operator. 
Let us briefly describe its construction. 

  Let $(M^{2l},\omega)$ be a symplectic manifold. Suppose there exists a metaplectic $(\Lambda, \mathcal{P})$ structure on $(M,\omega).$ (The definition of a metaplectic structure on a symplectic manifold is completely parallel to the Riemannian case. For this notion  see, e.g., Habermann, Habermann \cite{HH}.) 
Let us consider the associated vector bundle $\mathcal{S}$ to the principal $Mp(2l,\mathbb{R})$-bundle $\mathcal{P}$ via the metaplectic representation of the metaplectic group on the space of symplectic spinors, i.e.,
$\mathcal{S}:=\mathcal{P}\times_{meta} {\bf S}.$  This bundle is called symplectic spinor bundle. Let $\nabla$ be an affine symplectic torsion-free connection on $(M,\omega).$ 
Then the symplectic Dirac operator $\mathfrak{D}$ is uniquely defined as the composition of the symplectic Clifford multiplication and the associated covariant derivative on the symplectic spinor bundle.
The symplectic Rarita-Schwinger operator $\mathfrak{R}$  and other higher symplectic spin operators  could be defined by associating the studied $Mp(\mathbb{V},\omega)$-module ${\bf W}=\bigwedge^{\bullet}\mathbb{V}^*\otimes {\bf S}$ to the metaplectic structure. See Kr\'ysl, \cite{SK} for details.  
For a comprehensive introduction to the symplectic Dirac operators, see Habermann, Habermann \cite{HH}.
  
 Further, let us denote the spectrum of an endomorphism $\mathfrak{A}$ by $\hbox{Spec}(\mathfrak{A})$ and the space of eigenvectors of $\mathfrak{A}$ by $Eigen(\mathfrak{A})$.
 
The following notion is introduced for some technical reasons.
We call a non-zero section $\psi \in \Gamma(M,\mathcal{S})$ symplectic Killing spinor if there exists a complex number $\mu \in \mathbb{C}$ such that $\nabla_{X}\psi=\mu X.\psi$ for all $X \in \Gamma(M, TM).$ We call the complex number $\mu$ symplectic Killing number. Let us denote the set of symplectic Killing number by $\hbox{Kill}$ and the set of symplectic Killing spinors by $Kill.$
In the following theorem, the spectra of the symplectic Dirac and the symplectic Rarita-Schwinger are related.
The operator $T_1$ in the formulation of the theorem is the so called symplectic twistor operator.
 
{\bf Theorem 21:} Let $(M,\omega)$ be a symplectic manifold of dimension $2l$ admitting a metaplectic structure and $\nabla$ be a flat symplectic torsion-free connection. Then
\begin{itemize}
\item[1.)] If $\lambda \in \hbox{Spec}(\mathfrak{D})\setminus -\imath l (\hbox{Kill})$ then $\lambda -\frac{1}{l}\lambda \in \hbox{Spec}(\mathfrak{R}).$
\item[2.)] If $\psi \in Eigen(\mathfrak{D})\setminus Kill$ then
$T_1 \psi \in Eigen(\mathfrak{R}).$
\end{itemize} 
{\it Proof.} See Kr\'{y}sl \cite{SK}.

{\bf Remark:} In the proof, the relations between the basis elements $\{e^{\pm},h,f^{\pm}\}$ of the super Lie algebra
$\mathfrak{osp}(1|2)$ were used. The use of these relations makes it possible to almost avoid coordinate-computations.

The studied representation $\rho$ of the group $Mp(\mathbb{V},\omega)$ could be also used in a classification of invariant first order differential operators in the so called contact projective geometry (a type of Cartan geometry). See Kr\'{y}sl  \cite{Krysl} for details.

\end{document}